\documentclass{article}%
\usepackage{amssymb}
\usepackage{amsfonts}
\usepackage{amsmath}%
\usepackage{graphicx}

\begin{document}

\title{Marginalized Particle Filtering and Related Filtering Techniques as Message Passing}

\maketitle

\begin{abstract}
In this manuscript a factor graph approach is employed to investigate the
recursive filtering problem for mixed linear/nonlinear state-space models. Our
approach allows us to show that: a) the factor graph characterizing the
considered filtering problem is not cycle free; b) in the case of
conditionally linear Gaussian systems, applying the sum-product rule, together
with different scheduling procedures for message passing, to this graph
results in both known and novel filtering techniques. In particular, it is
proved that, on the one hand, adopting a specific message scheduling for
forward only message passing leads to marginalized particle filtering in a
natural fashion; on the other hand, if iterative strategies for message
passing are employed, a novel filtering method, dubbed turbo filter for its
conceptual resemblance to the turbo decoding methods devised for concatenated
channel codes, can be developed.

\end{abstract}

\bigskip

\begin{center}
Giorgio M. Vitetta, Emilio Sirignano, Francesco Montorsi and Matteo Sola

\vspace{5mm}University of Modena and Reggio Emilia

Department of Engineering "Enzo Ferrari"

Via P. Vivarelli 10/1, 41125 Modena - Italy

email: giorgio.vitetta@unimore.it, emilio.sirignano@unimore.it,
francesco.montorsi@gmail.com, sola.matteo87@gmail.com
\end{center}

\bigskip


\pagestyle{myheadings} \markright{G. M. Vitetta et al.,  Marginalized Particle
Filtering \ldots}

\bigskip\vspace{1cm}

\baselineskip0.2 in\newpage

\section{Introduction\label{sec:intro}}

The nonlinear filtering problem consists of inferring the posterior
distribution of the hidden state of a nonlinear dynamic system from a set of
past and present measurements \cite{Arulampalam_2002}. It is well known that,
if a nonlinear dynamic system can be described by a \emph{state-space model}
(SSM), a general sequential procedure, based on the Bayes' rule and known as
\emph{Bayesian filtering}, can be easily derived for recursively computing the
posterior distribution of the system current state \cite{Arulampalam_2002}.
Unluckily, Bayesian filtering is analytically tractable in few cases for the
following two reasons \cite{Daum_1986}: a) one of the two steps it consists of
requires multidimensional integration which, in most cases, does not admit a
closed form solution; b) the functional form of the required posterior
distribution may not be preserved over successive recursions. For this reason,
sequential techniques employed in practice are based on various analytical
approximations and, consequently, generate a functional approximation of the
desired distribution. Such techniques are commonly divided into \emph{local}
and \emph{global} methods on the basis of the way posterior distributions are
approximated \cite{Mazuelas_2013,Smidl_2008,Simandla_2006}. Local methods,
like \emph{extended Kalman filtering} \cite{Anderson_1979} and \emph{unscented
filtering} \cite{Julier_2004}, are computationally efficient, but may suffer
from the problem of error accumulation over time. On the contrary, global
methods, like \emph{sequential Monte Carlo methods}
\cite{Doucet_2001,Doucet_2000} (also known as \emph{particle filtering}, PF,
methods \cite{Gustafsson_2010,Gustafsson_et_al_2002,Nordlund_2009}) and
\emph{point mass filtering} \cite{Simandla_2006,Bucy_1971} may achieve high
accuracy at the price, however, of an unmanageable complexity and numerical
problems in the presence of a large dimension of system state \cite{Daum_2003}%
. These considerations have motivated various research activities focused on
the development of novel Bayesian filters able to achieve high accuracy under
given computational constraints. Significant results in this research area
concern the use of the new representations for complex distributions, like
\emph{belief condensation filtering} \cite{Mazuelas_2013}, and the development
of novel filtering techniques combining local and global methods, like
\emph{marginalized particle filtering} (MPF)
\cite{Schon_2005,Schon_2005_complexity}, and other methods originating from it
\cite{Smidl_2008,Mustiere_2006,Lu_2007}. Note that the last class of methods
applies to \emph{mixed nonlinear/nonlinear models} \cite{Lindsten_2016}, that
is to models whose state can be partitioned in a conditionally linear portion
(usually called \emph{linear state variable}) and in a nonlinear portion
(representing the remaining part of system state and called \emph{nonlinear
state variable}). This partitioning of system state allows to combine a global
method (e.g., particle filtering) operating on the nonlinear state variable
with a local technique (e.g., Kalman filtering) involving the linear state
variable only.

In this manuscript the \emph{factor graph} (FG) approach illustrated by
Loeliger \emph{et al}. in \cite{Loeliger_2007} is employed to revisit the
problem of recursive Bayesian filtering for mixed linear/nonlinear models from
a perspective substantially different from that adopted in MPF
\cite{Schon_2005}. This allows us to shed new light on the problem of
filtering for mixed linear/nonlinear models, providing a new interpretation of
MPF and paving the way for the development of new filtering techniques. In
particular, based on this approach, we are able to show that: a) the
considered filtering problem can be formulated as a message massing problem
over a specific FG, which, unluckily, is \emph{not cycle free}; b) in the case
of a \emph{conditionally linear Gaussian} (CLG) SSM \cite{Lindsten_2016}, MPF
results from the application of the \emph{sum-product algorithm} (SPA)
\cite{Loeliger_2007,Kschischang_2001}, together with a specific scheduling
procedure for forward only message passing, to this graph; c) our graphical
representation leads, in a natural fashion, to the development of novel
filtering methods simplifying and/or generalising it. As far as the last point
is concerned, in our work specific attention is paid to the development of a
novel \emph{iterative filtering technique} that exploits the exchange of
probabilistic (i.e., soft) messages to progressively refine the posteriors of
the linear and nonlinear state variables within each recursion and is dubbed
\emph{turbo filtering} (TF) for its conceptual resemblance to the iterative
(i.e., turbo) decoding of concatenated channel codes.

It is important to point that our approach has been inspired by various ideas
and results already available in the technical literature concerning different
research areas; here, we limit to mention the following relevant facts:

\begin{itemize}
\item A mixed linear/nonlinear Markov system can be represented as the
\emph{concatenation} of two interacting subsystems, one governed by linear
dynamics, the other one accounting for a nonlinear behavior; conceptually
related (finite state) Markov models can be found in data communications and,
in particular, in concatenated channel coding (e.g., turbo coding
\cite{Berrou_1996}) and in coded transmissions over inter-symbol interference
channels for which \emph{turbo decoding methods}
\cite{Berrou_1996,Benedetto_1998} and \emph{turbo equalization techniques}
\cite{Koetter_2004} have been developed, respectively\footnote{Note that these
classes of algorithms can be seen as specific applications of the so called
\emph{turbo principle} \cite{Hagenauer_1997}, \cite[Par. 10.5.1]{Vitetta}}.

\item Factor graphs play an essential role in the derivation and
interpretation of turbo decoding and equalization
\cite{Loeliger_2007,Worthen_2001} (for instance, turbo decoding techniques
emerge in a natural fashion from graphical models of codes
\cite{Kschischang_1998}).

\item Both Kalman filtering and particle methods can be viewed as
\emph{message passing procedures on factor graphs}, as shown in
\cite{Loeliger_2007,Kschischang_2001} and in \cite{Dauwels_2006}, respectively.

\item Various methods to progressively refine distributional approximations
through multiple iterations have been developed in the field of Bayesian
inference on dynamic systems (even if implementations substantially different
from that we devise have been proposed), and, in particular,
in\emph{\ expectation propagation} in Bayesian networks
\cite{Minka_2001,Zoeter_2006} and in \emph{variational Bayesian filtering}
\cite{Smidl_2008}. Consequently, various links to previous work on
\emph{Bayesian inference on graphical models} and \emph{variational Bayes
methods} \cite{Smidl_2005} can be also established.
\end{itemize}

The remaining part of this manuscript is organized as follows. The
mathematical model of the considered class of mixed linear/nonlinear systems
is illustrated in Section \ref{sec:scenario}, whereas a representation of the
filtering problem for these systems through a proper FG is provided in Section
\ref{sec:Factorgraphs}. Then, it is shown that applying the SPA and proper
message scheduling strategies to this FG leads to MPF in Section
\ref{sec:Message-Passing}. This approach paves the way, in a natural fashion,
for the development of simplifications and generalizations of MPF, that are
devised in Sections \ref{sec:simplifying} and
\ref{sec:Modifications-and-Extensions}, respectively. The novel filtering
methods proposed in this manuscript are compared, in terms of accuracy and
computational effort, with MPF in Section \ref{num_results}. Finally, some
conclusions are offered in Section \ref{sec:conc}.

\emph{Notations}: The \emph{probability density function} (pdf) of a random
vector $\mathbf{R}$ evaluated at point $\mathbf{r}$ is denoted $f(\mathbf{r})
$; $\mathcal{N}\left(  \mathbf{r};\mathbf{\eta_{r}},\mathbf{C_{r}}\right)  $
represents the pdf of a Gaussian random vector $\mathbf{R}$ characterized by
the mean $\mathbf{\eta_{r}}$ and covariance matrix $\mathbf{\mathbf{C_{r}}}$
evaluated at point $\mathbf{r}$; the \emph{precision }(or \emph{weight})
\emph{matrix} associated with the covariance matrix $\mathbf{\mathbf{C_{r}}}$
is denoted $\mathbf{\mathbf{W_{r}}}$, whereas the \emph{transformed mean
vector} $\mathbf{\mathbf{W_{r}}\eta_{r}}$ is denoted $\mathbf{\mathbf{w_{r}}}%
$; $x_{i}$ denotes the $i$-th element of the vector $\mathbf{x}$.

\section{System Model\label{sec:scenario}}

In the following we focus on a discrete-time mixed linear/nonlinear SSM
\cite{Schon_2005}, whose \emph{hidden state} in the $l$-th interval is
represented by a $D$-dimensional real vector $\mathbf{x}_{l}\triangleq\lbrack
x_{0,l},x_{1,l},...,$ $x_{D-1,l}]^{T}$. We assume that this vector can be
partitioned as
\begin{equation}
\mathbf{x}_{l}=\left[  \left(  \mathbf{x}_{l}^{(L)}\right)  ^{T},\left(
\mathbf{x}_{l}^{(N)}\right)  ^{T}\right]  ^{T},\label{eq:state_structure}%
\end{equation}
where $\mathbf{x}_{l}^{(L)}\triangleq\lbrack x_{0,l}^{(L)},x_{1,l}%
^{(L)},...,x_{D_{L}-1,l}^{(L)}]^{T}$ ($\mathbf{x}_{l}^{(N)}\triangleq\lbrack
x_{0,l}^{(N)},x_{1,l}^{(N)},...,x_{D_{N}-1,l}^{(L)}]^{T}$) is the so called
\emph{linear }(\emph{nonlinear}) \emph{component} of $\mathbf{x}_{l}$
(\ref{eq:state_structure}), with $D_{L}<D$ ($D_{N}=D-D_{L}$). This
partitioning of $\mathbf{x}_{l}$ is accomplished as follows. First,
$\mathbf{x}_{l}^{(L)}$ is identified as that portion of $\mathbf{x}_{l}$
(\ref{eq:state_structure}) characterized by the following two properties:

\begin{enumerate}
\item \emph{Conditionally linear dynamics} - This means that its update
equation, conditioned on $\mathbf{x}_{l}^{(N)}$, is \emph{linear}, so that
\begin{equation}
\mathbf{x}_{l+1}^{(L)}=\mathbf{A}_{l}^{(L)}\left(  \mathbf{x}_{l}%
^{(N)}\right)  \mathbf{x}_{l}^{(L)}+\mathbf{f}_{l}^{(L)}\left(  \mathbf{x}%
_{l}^{(N)}\right)  +\mathbf{w}_{l}^{(L)},\label{eq:XL_update}%
\end{equation}
where $\mathbf{f}_{l}^{(L)}\left(  \mathbf{x}\right)  $ is a time-varying
$D_{L}$-dimensional real function, $\mathbf{A}_{l}^{(L)}(\mathbf{x}_{l}^{(N)})
$ is a time-varying $D_{L}\times D_{L}$ real matrix and $\mathbf{w}_{l}^{(L)}
$ is the $l$-th element of the process noise sequence $\{\mathbf{w}_{k}%
^{(L)}\}$, which consists of $D_{L}$- dimensional \emph{independent and
identically distributed} (iid) noise\emph{\ }vectors.

\item \emph{Conditionally linear} (or \emph{almost linear}) dependence of all
the available measurements on it - In other words, these quantities,
conditioned on $\mathbf{x}_{l}^{(N)}$, exhibit a \emph{linear} dependence on
$\mathbf{x}_{l}^{(L)}$ (additional details about this feature are provided below).
\end{enumerate}

Then, $\mathbf{x}_{l}^{(N)}$ is generated by putting together all the
components of $\mathbf{x}_{l}$ that do not belong to $\mathbf{x}_{l}^{(L)}$.
For this reason, generally speaking, this vector is characterized by at least
one of the following two properties:

a) \emph{Nonlinear dynamics} - The update equation
\begin{equation}
\mathbf{x}_{l+1}^{(N)}=\mathbf{f}_{l}^{(N)}\left(  \mathbf{x}_{l}%
^{(N)}\right)  +\mathbf{A}_{l}^{(N)}\left(  \mathbf{x}_{l}^{(N)}\right)
\mathbf{x}_{l}^{(L)}+\mathbf{w}_{l}^{(N)}\label{eq:XN_update}%
\end{equation}
is assumed in the following for the nonlinear component of system state, where
$\mathbf{A}_{l}^{(N)}(\mathbf{x}_{l}^{(N)})$ is a time-varying $D_{N}\times
D_{L}$ real matrix, $\mathbf{f}_{l}^{(N)}\left(  \mathbf{x}\right)  $ is a
time-varying $D_{N}$-dimensional real function and $\mathbf{w}_{l}^{(N)}$ is
the $l$-th element of the process noise sequence $\{\mathbf{w}_{k}^{(N)}\}$,
which consists of $D_{N}$-dimensional iid noise vectors and is statistically
independent of $\{\mathbf{w}_{k}^{(L)}\}$.

b) \emph{A nonlinear} \emph{dependence of all the available measurements on
it} (further details are provided below).

In the following Section we focus on the so-called \emph{filtering problem},
which concerns the evaluation of the posterior pdf $f(\mathbf{x}%
_{t}|\mathbf{y}_{1:t})$ at an instant $t>1$, given a) the initial pdf
$f(\mathbf{x}_{1})$ and b) the $t\cdot P$-dimensional \emph{measurement}
vector
\begin{equation}
\mathbf{y}_{1:t}=\left[  \mathbf{y}_{1}^{T},\mathbf{y}_{2}^{T},...,\mathbf{y}%
_{t}^{T}\right]  ^{T},\label{eq:y_0:N}%
\end{equation}
where $\mathbf{y}_{l}\triangleq\lbrack y_{0,l},y_{1,l},$ $...,y_{P-1,l}]^{T} $
denotes the $P$-dimensional real vector collecting all the noisy measurements
available at time $l$. As already mentioned above, the measurement vector
$\mathbf{y}_{l}$ exhibits a \emph{linear} (\emph{nonlinear}) dependence on
$\mathbf{x}_{l}^{(L)}$ ($\mathbf{x}_{l}^{(N)}$), so that the model
\cite{Lu_2007}
\begin{equation}
\mathbf{y}_{l}=\mathbf{h}_{l}\left(  \mathbf{x}_{l}^{(N)}\right)
+\mathbf{B}_{l}\left(  \mathbf{x}_{l}^{(N)}\right)  \mathbf{x}_{l}%
^{(L)}+\mathbf{e}_{l}\label{eq:y_t}%
\end{equation}
can be adopted, where $\mathbf{B}_{l}(\mathbf{x}_{l}^{(N)})$ is a time-varying
$P\times D_{L}$ real matrix, $\mathbf{h}_{l}(\mathbf{x}_{l}^{(N)})$ is a
time-varying $P$-dimensional real function and $\mathbf{e}_{l}$ the $l$-th
element of the measurement noise sequence $\left\{  \mathbf{e}_{k}\right\}  $
consisting of $P$-dimensional iid noise vectors and independent of both
$\{\mathbf{w}_{k}^{(N)}\}$ and $\{\mathbf{w}_{k}^{(L)}\}$.

\section{Representation of the Filtering Problem via Factor
Graphs\label{sec:Factorgraphs}}

\begin{figure}
	\centering
		\includegraphics[width=0.55\textwidth]{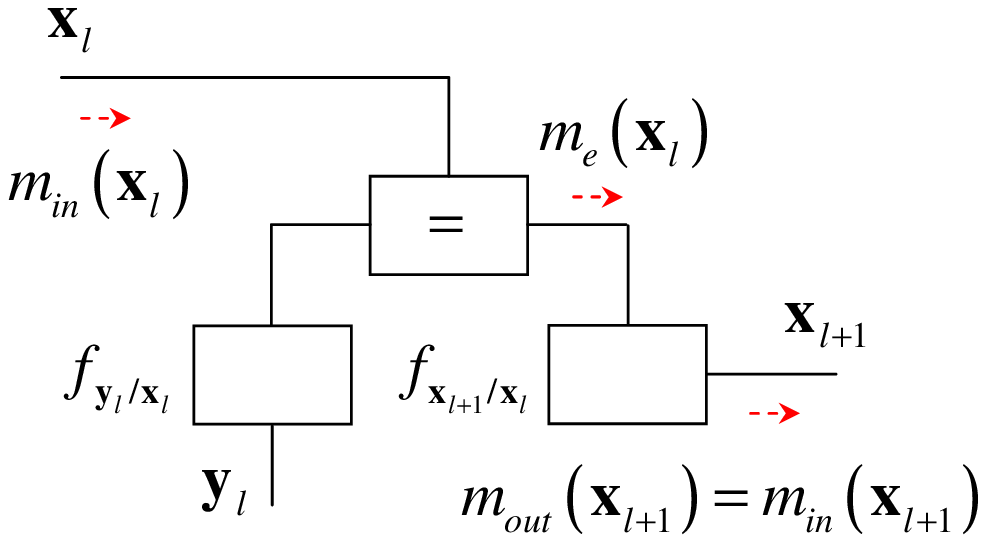}
	\caption{Factor graph representing (\ref{eq:meas_update_m}) and
(\ref{eq:time_update_m}). The SPA message flow characterizing the $l$-th
recursion of Bayesian filtering is indicated by red arrows.}
	\label{Fig_1}
\end{figure}

Generally speaking, the \emph{filtering problem} for a SSM described by the
\emph{Markov model} $f(\mathbf{x}_{l+1}|\mathbf{x}_{l})$ and the
\emph{observation model} $f(\mathbf{y}_{l}|\mathbf{x}_{l})$ for any $l$
concerns the computation of the posterior pdf $f(\mathbf{x}_{t}|\mathbf{y}%
_{1:t})$ for $t\geq1$ by means of a recursive procedure
\cite{Arulampalam_2002}. It is well known that, if the pdf $f(\mathbf{x}_{1})$
is known, a general \emph{Bayesian recursive procedure}, consisting of a
\emph{measurement update} step followed by a \emph{time update} step, can be
employed. In practice, in the first step of the $l$-th recursion (with
$l=1,2,...,t$) the conditional pdf
\begin{equation}
f\left(  \mathbf{x}_{l}\left\vert \mathbf{y}_{1:l}\right.  \right)  =f\left(
\mathbf{x}_{l}\left\vert \mathbf{y}_{1:(l-1)}\right.  \right)  f\left(
\mathbf{y}_{l}\left\vert \mathbf{x}_{l}\right.  \right)  \frac{1}{f\left(
\mathbf{y}_{l}\left\vert \mathbf{y}_{1:(l-1)}\right.  \right)  }%
\label{eq:meas_update}%
\end{equation}
is computed on the basis of pdf $f(\mathbf{x}_{l}|\mathbf{y}_{1:(l-1)})$
(evaluated in the last step of the previous recursion\footnote{Note that in
the first recursion (i.e., for $l=1$) $f(\mathbf{x}_{l}|\mathbf{y}%
_{1:(l-1)})=f(\mathbf{x}_{1}|\mathbf{y}_{1:0})=f(\mathbf{x}_{1})$ and
$f(\mathbf{y}_{l}|\mathbf{y}_{1:(l-1)})=f(\mathbf{y}_{1}|\mathbf{y}%
_{1:0})=f(\mathbf{y}_{1})$, so that $f(\mathbf{x}_{1}|\mathbf{y}%
_{1})=f(\mathbf{x}_{1})f(\mathbf{y}_{1}|\mathbf{x}_{1})/f(\mathbf{y}_{1})$.}),
the present measurement vector $\mathbf{y}_{l}$ and the pdf
\begin{equation}
f\left(  \mathbf{y}_{l}\left\vert \mathbf{y}_{1:(l-1)}\right.  \right)  =\int
f\left(  \mathbf{y}_{l}\left\vert \mathbf{x}_{l}\right.  \right)  f\left(
\mathbf{x}_{l}\left\vert \mathbf{y}_{1:(l-1)}\right.  \right)  d\mathbf{x}%
_{l}.\label{eq:norm_factor}%
\end{equation}
In the second step $f(\mathbf{x}_{l}|\mathbf{y}_{1:l})$ (\ref{eq:meas_update})
is exploited to compute the pdf
\begin{equation}
f\left(  \mathbf{x}_{l+1}\left\vert \mathbf{y}_{1:l}\right.  \right)  =\int
f\left(  \mathbf{x}_{l+1}\left\vert \mathbf{x}_{l}\right.  \right)  f\left(
\mathbf{x}_{l}\left\vert \mathbf{y}_{1:l}\right.  \right)  d\mathbf{x}%
_{l},\label{eq:time_update}%
\end{equation}
which represents a \emph{prediction} about the future state $\mathbf{x}_{l+1}
$. It is important to point out that: 1) the term $1/f(\mathbf{y}%
_{l}|\mathbf{y}_{1:(l-1)})$ appearing in the \emph{right hand side} (RHS) of
(\ref{eq:meas_update}) represents a \emph{normalization factor}; 2) both
(\ref{eq:norm_factor}) and (\ref{eq:time_update}) require integration with
respect to $\mathbf{x}_{l}$ and this may represent a formidable task when the
dimensionality of $\mathbf{x}_{l}$ is large and/or the pdfs appearing in the
integrands are not Gaussian; 3) this recursive procedure lends itself to be
efficiently represented by a \emph{message passing algorithm} over a
proper\emph{\ }FG\footnote{\emph{Forney-style} factor graphs are always
considered in the following \cite{Loeliger_2007}.} \cite{Dauwels_2006}, in
which each \emph{factor} of a product of functions is represented by a
distinct \emph{node} (a rectangle in our diagrams), whereas each
\emph{variable} is associated with a specific (and, usually, unoriented)
\emph{edge} or \emph{half edge}. As far as the last point is concerned, we
also note that the derivation of this FG relies on the fact that the a
posteriori pdf $f(\mathbf{x}_{t}|\mathbf{y}_{1:t})$ has the \emph{same }FG as
the joint pdf $f(\mathbf{x}_{t},\mathbf{y}_{1:t})$ (see \cite[Sec. II, p.
1297]{Loeliger_2007}) and the last pdf can be computed recursively through a
procedure similar to that illustrated above, but in which the measurement
update (\ref{eq:meas_update}) and the time update (\ref{eq:time_update}) are
replaced by
\begin{equation}
f\left(  \mathbf{x}_{l},\mathbf{y}_{1:l}\right)  =f\left(  \mathbf{x}%
_{l},\mathbf{y}_{1:(l-1)}\right)  f\left(  \mathbf{y}_{l}\left\vert
\mathbf{x}_{l}\right.  \right)  ,\label{eq:meas_update_m}%
\end{equation}
and
\begin{equation}
f\left(  \mathbf{x}_{l+1},\mathbf{y}_{1:l}\right)  =\int f\left(
\mathbf{x}_{l+1}\left\vert \mathbf{x}_{l}\right.  \right)  f\left(
\mathbf{x}_{l},\mathbf{y}_{1:l}\right)  d\mathbf{x}_{l}%
,\label{eq:time_update_m}%
\end{equation}
respectively, so that the evaluation of the above mentioned normalization
factor is no more required. In fact, eqs. (\ref{eq:meas_update_m}) and
(\ref{eq:time_update_m}) involve only \emph{products of pdfs} and a \emph{sum}
(i.e., integration) \emph{of products}, so that they can be represented by
means of the FG shown in Fig. \ref{Fig_1} (where, following
\cite{Loeliger_2007}, a simplified notation is used for the involved pdfs and
the \emph{equality constraint node}, which represents an equality constraint
\textquotedblleft function\textquotedblright, that is a Dirac delta function).
Since this FG is \emph{cycle free}, the pdf $f(\mathbf{x}_{l},\mathbf{y}%
_{1:l})$ can be evaluated applying the well known \emph{sum-product
rule}\footnote{In a Forney-style FG, such a rule can be formulated as follows
\cite{Loeliger_2007}: the message emerging from a node \emph{f} along some
edge \emph{x} is formed as the product of \emph{f} and all the incoming
messages along all the edges that enter the node \emph{f} except \emph{x},
summed over all the involved variables except \emph{x}.} (i.e., the SPA) to
it, i.e. developing a proper mechanism for passing probabilistic messages
along this FG (the flow of messages is indicated by red arrows in Fig.
\ref{Fig_1}). In fact, if the input message $\vec{m}_{in}\left(
\mathbf{x}_{l}\right)  =f(\mathbf{x}_{l},\mathbf{y}_{1:(l-1)})$ enters this
FG, the message going out of the \emph{equality node} is given by
\begin{equation}
\vec{m}_{e}\left(  \mathbf{x}_{l}\right)  =\vec{m}_{in}\left(  \mathbf{x}%
_{l}\right)  f\left(  \mathbf{y}_{l}\left\vert \mathbf{x}_{l}\right.  \right)
,\label{eq:message_y}%
\end{equation}
so that $\vec{m}_{e}\left(  \mathbf{x}_{l}\right)  =f(\mathbf{x}%
_{l},\mathbf{y}_{1:l})$ (see (\ref{eq:meas_update_m})); then, the message
emerging from the \emph{function node} referring to the pdf $f(\mathbf{x}%
_{l+1}|\mathbf{x}_{l}) $ is expressed by
\begin{equation}
\vec{m}_{out}\left(  \mathbf{x}_{l+1}\right)  =\int f\left(  \mathbf{x}%
_{l+1}\left\vert \mathbf{x}_{l}\right.  \right)  \vec{m}_{e}\left(
\mathbf{x}_{l}\right)  d\mathbf{x}_{l},\label{eq:message_y-1}%
\end{equation}
so that $\vec{m}_{out}\left(  \mathbf{x}_{l+1}\right)  =f(\mathbf{x}%
_{l+1},\mathbf{y}_{1:l})=\vec{m}_{in}\left(  \mathbf{x}_{l+1}\right)  $ (see
(\ref{eq:time_update_m})). From this result it can be easily inferred that the
pdf $f(\mathbf{x}_{t},\mathbf{y}_{1:t})$ (and, up to a scale factor, the pdf
$f(\mathbf{x}_{t}|\mathbf{y}_{1:t})$) results from the application of the SPA
to the overall FG originating from the ordered concatenation of multiple
subgraphs, each structured like the one shown in Fig. \ref{Fig_1} and
associated with $l=1,2,...,t$. In this graph the flow of messages produced by
the SPA proceeds from left to right, i.e. the pdf $f(\mathbf{x}_{t}%
,\mathbf{y}_{1:t})$ is generated by a \emph{forward only} message passing.
Note also that, in principle, the desired pdf $f(\mathbf{x}_{t},\mathbf{y}%
_{1:t})$ is computed as the product between two messages, one for each
direction, reaching the rightmost \emph{half edge} of the overall FG, but one
of the two incoming messages for that edge is the constant function
$\overleftarrow{m}_{he}(\mathbf{x}_{t})=1$.%

Unluckily, as the size $D$ of $\mathbf{x}_{l}$ (\ref{eq:state_structure}) gets
large, the computational burden associated with (\ref{eq:meas_update}%
)-(\ref{eq:time_update}) (or, equivalently, (\ref{eq:meas_update_m}) and
(\ref{eq:time_update_m})) becomes unmanageable. In principle, a substantial
complexity reduction can be achieved \emph{decoupling}\footnote{Note that,
generally speaking, in the considered problem the coupling of the filtering
problem for $\mathbf{x}_{l}^{(L)}$ with that for $\mathbf{x}_{l}^{(N)}$ is due
not only to the structure of the update equations (\ref{eq:XL_update}) and
(\ref{eq:XN_update}), but also to the measurement vector $\mathbf{y}_{l}$
(\ref{eq:y_t}), since this exhibits a mixed dependence on the two components
of the state vector $\mathbf{x}_{l}$ (\ref{eq:state_structure}).} the
filtering problem for $\mathbf{x}_{l}^{(L)}$ from that for $\mathbf{x}%
_{l}^{(N)}$, i.e. the evaluation of $f(\mathbf{x}_{l}^{(L)}|\mathbf{y}_{1:l})$
from that of $f(\mathbf{x}_{l}^{(N)}|\mathbf{y}_{1:l})$. In fact, this
approach potentially provides the following two benefits: a) a given filtering
problem is turned into a couple of filtering problems of smaller
dimensionality and b) some form of computationally efficient standard
filtering (e.g., Kalman or extended Kalman filtering) can be hopefully
exploited for the linear portion $\mathbf{x}_{l}^{(L)}$ of the state vector
$\mathbf{x}_{l}$ (\ref{eq:state_structure}). As a matter of fact, these
principles have been exploited in devising MPF
\cite{Gustafsson_2010,Schon_2005} and, as it will become clearer in the
following, they must be always kept into account in the derivation of our FG
representation. Before illustrating this derivation, however, the measurement
and state models on which such a representation relies need to be clearly
defined; for this reason, these models are analysed in detail in the following
part of this Section. To begin, let us concentrate on the models involved in
the filtering problem for $\mathbf{x}_{l}^{(L)}$, i.e. on the evaluation of
$f(\mathbf{x}_{l}^{(L)}|\mathbf{y}_{1:l})$, under the assumption that the
nonlinear portion $\mathbf{x}_{l}^{(N)}$ of the system state is \emph{known}
for any $l$. In this case, the evaluation of this pdf can benefit not only
from the knowledge of $\mathbf{y}_{l}$ (\ref{eq:y_t}), but also from that of
the quantity (see (\ref{eq:XN_update}))%
\begin{equation}
\mathbf{z}_{l}^{(L)}\triangleq\mathbf{x}_{l+1}^{(N)}-\mathbf{f}_{l}%
^{(N)}\left(  \mathbf{x}_{l}^{(N)}\right)  =\mathbf{A}_{l}^{(N)}\left(
\mathbf{x}_{l}^{(N)}\right)  \mathbf{x}_{l}^{(L)}+\mathbf{w}_{l}%
^{(N)},\label{eq:z_L_l}%
\end{equation}
which can be interpreted as a \emph{pseudo-measurement} \cite{Schon_2005},
since it does not originate from real measurements, but from the constraints
expressed by the state equation (\ref{eq:XN_update}). This leads to
considering the \emph{overall observation model}
\begin{equation}%
\begin{array}
[c]{c}%
f\left(  \mathbf{y}_{l},\mathbf{z}_{l}^{(L)}\left\vert \mathbf{x}_{l}%
^{(L)},\mathbf{x}_{l}^{(N)}\right.  \right) \\
=f\left(  \mathbf{y}_{l}\left\vert \mathbf{x}_{l}^{(L)},\mathbf{x}_{l}%
^{(N)}\right.  \right)  f\left(  \mathbf{z}_{l}^{(L)}\left\vert \mathbf{x}%
_{l}^{(L)},\mathbf{x}_{l}^{(N)}\right.  \right)
\end{array}
\label{eq:meas_overall_L}%
\end{equation}
for $\mathbf{x}_{l}^{(L)}$, where
\begin{equation}
f\left(  \mathbf{y}_{l}\left\vert \mathbf{x}_{l}^{(L)},\mathbf{x}_{l}%
^{(N)}\right.  \right)  =\left.  f\left(  \mathbf{e}_{l}\right)  \right\vert
_{\mathbf{e}_{l}=\mathbf{y}_{l}-\mathbf{B}_{l}\left(  \mathbf{x}_{l}%
^{(N)}\right)  \mathbf{x}_{l}^{(L)}-\mathbf{h}_{l}\left(  \mathbf{x}_{l}%
^{(N)}\right)  },\label{eq:f_y_L_cond_x_L-1}%
\end{equation}
and
\begin{equation}
f\left(  \mathbf{z}_{l}^{(L)}\left\vert \mathbf{x}_{l}^{(L)},\mathbf{x}%
_{l}^{(N)}\right.  \right)  =\left.  f\left(  \mathbf{w}_{l}^{(N)}\right)
\right\vert _{\mathbf{w}_{l}^{(N)}=\mathbf{z}_{l}^{(L)}-\mathbf{A}_{l}%
^{(N)}\left(  \mathbf{x}_{l}^{(N)}\right)  \mathbf{x}_{l}^{(L)}}%
.\label{eq:f_z_L_cond_x_L-1}%
\end{equation}
If the observation model (\ref{eq:meas_overall_L}) and the \emph{state model}
(see (\ref{eq:XL_update}))
\begin{equation}%
\begin{array}
[c]{c}%
f\left(  \mathbf{x}_{l+1}^{(L)}\left\vert \mathbf{x}_{l}^{(L)},\mathbf{x}%
_{l}^{(N)}\right.  \right) \\
=f_{\mathbf{w}^{(L)}}\left(  \mathbf{x}_{l+1}^{(L)}-\mathbf{f}_{l}%
^{(L)}\left(  \mathbf{x}_{l}^{(N)}\right)  -\mathbf{A}_{l}^{(L)}\left(
\mathbf{x}_{l}^{(N)}\right)  \mathbf{x}_{l}^{(L)}\right)
\end{array}
\label{eq:state_model_L-1}%
\end{equation}
are adopted for $\mathbf{x}_{l}^{(L)}$, the graph identified by the blue lines
and rectangles appearing in Fig. \ref{Fig_1} can be drawn. Then, in principle,
if the sum-product rule is applied to it under the assumption that the couple
$(\mathbf{x}_{l}^{(N)},\mathbf{x}_{l+1}^{(N)})$ is known for any $l$, the
expressions of the messages flowing in the overall graph for the evaluation of
$f(\mathbf{x}_{t}^{(L)},\mathbf{y}_{1:t},\mathbf{z}_{1:t}^{(L)})$ can be
easily derived. It is important to point out that:

\begin{itemize}
\item This graph contains a node which does not refer to the above mentioned
density factorizations\footnote{This peculiarity is also evidenced by the
presence of an arrow on all the edges connected to such a node.}, but
represents the transformation from the couple $(\mathbf{x}_{l}^{(N)}%
,\mathbf{x}_{l+1}^{(N)})$ to $\mathbf{z}_{l}^{(L)}$ (see (\ref{eq:z_L_l}));
this feature of the graph has to be carefully kept into account when deriving
message passing algorithms.

\item Generally speaking, the evaluation of the conditional pdf $f(\mathbf{z}%
_{l}^{(L)}|\mathbf{x}_{l}^{(L)},\mathbf{x}_{l}^{(N)})$ requires the knowledge
of the joint pdf of $\mathbf{x}_{l+1}^{(N)}$ and $\mathbf{x}_{l}^{(N)}$
conditioned on $\mathbf{x}_{l}^{(L)}$ (see (\ref{eq:z_L_l})).
\end{itemize}

\begin{figure}
	\centering
		\includegraphics[width=0.75\textwidth]{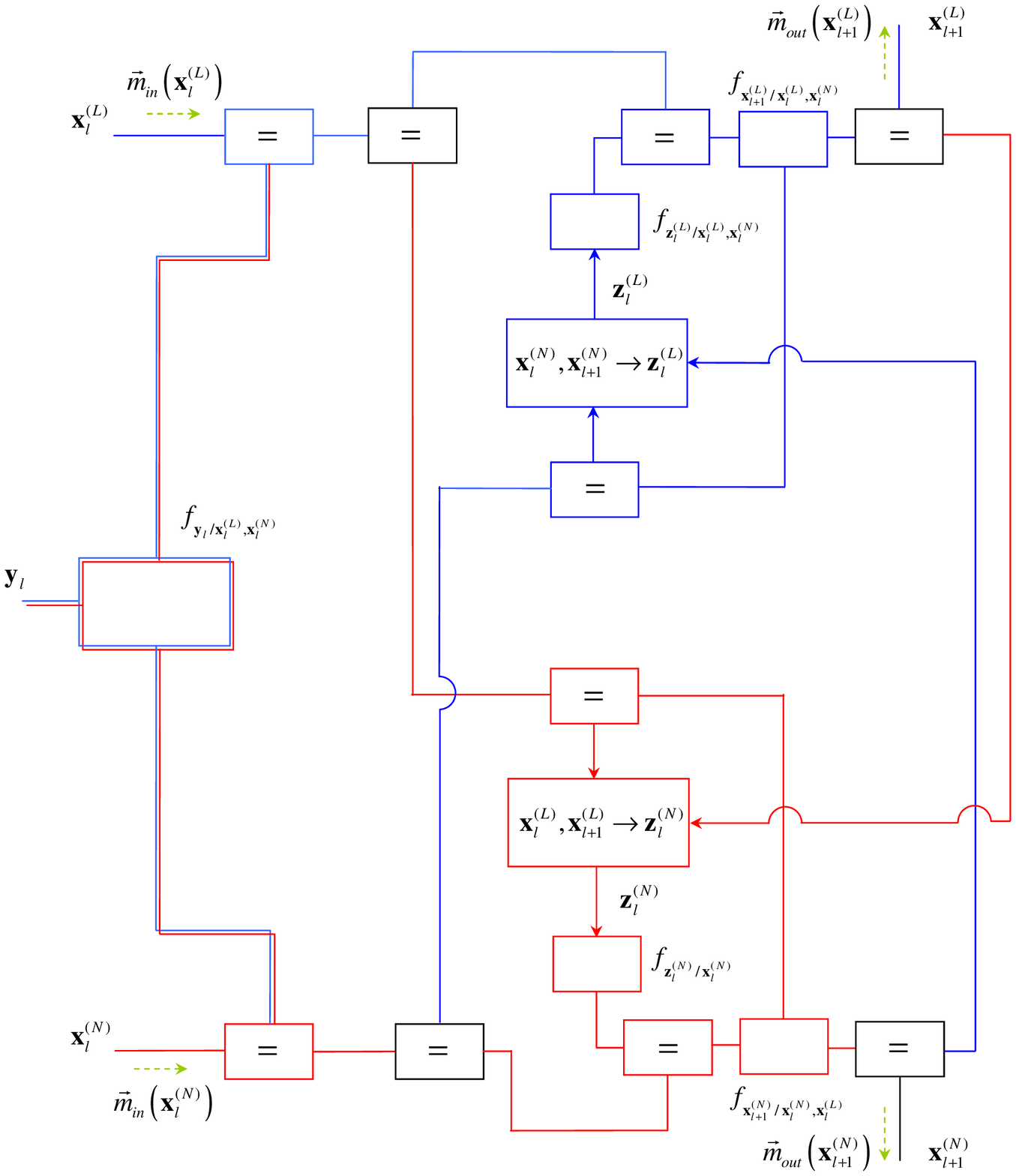}
	\caption{Overall factor graph resulting from the merge of two subgraphs, one
referring to filtering for $\mathbf{x}_{l}^{(L)}$ (in blue), the other one to
that for $\mathbf{x}_{l}^{(N)}$ (in red). \ The equality constraint nodes
introduced to connect these subgraphs are identified by black lines. The flow
of the messages along the half edges $\mathbf{x}_{l}^{(L)}$ and $\mathbf{x}%
_{l}^{(N)}$ (input) and that of the messages along the half edges
$\mathbf{x}_{l+1}^{(L)}$ and $\mathbf{x}_{l+1}^{(N)}$ (output) are indicated
by green arrows.}
	\label{Fig_2}
\end{figure}

The same line of reasoning can be followed for the filtering problem
concerning $\mathbf{x}_{l}^{(N)}$. Consequently, in this case the linear
portion $\mathbf{x}_{l}^{(L)}$ of the system state is assumed to be known for
any $l$ and the pseudo-measurement (see (\ref{eq:XL_update}))%
\begin{equation}
\mathbf{z}_{l}^{(N)}\triangleq\mathbf{x}_{l+1}^{(L)}-\mathbf{A}_{l}%
^{(L)}\left(  \mathbf{x}_{l}^{(N)}\right)  \mathbf{x}_{l}^{(L)}=\mathbf{f}%
_{l}^{(L)}\left(  \mathbf{x}_{l}^{(N)}\right)  +\mathbf{w}_{l}^{(L)}%
\label{eq:z_N_l}%
\end{equation}
is defined. This leads to the \emph{overall observation model}
\begin{equation}%
\begin{array}
[t]{c}%
f\left(  \mathbf{y}_{l},\mathbf{z}_{l}^{(N)}\left\vert \mathbf{x}_{l}%
^{(N)},\mathbf{x}_{l}^{(L)}\right.  \right)  =\\
f\left(  \mathbf{y}_{l}\left\vert \mathbf{x}_{l}^{(N)},\mathbf{x}_{l}%
^{(L)}\right.  \right)  f\left(  \mathbf{z}_{l}^{(N)}\left\vert \mathbf{x}%
_{l}^{(N)}\right.  \right)
\end{array}
\label{eq:meas_overall_N}%
\end{equation}
for $\mathbf{x}_{l}^{(N)}$, where $f(\mathbf{z}_{l}^{(N)}|\mathbf{x}_{l}%
^{(N)})$ can be expressed similarly as $f(\mathbf{z}_{l}^{(L)}|\mathbf{x}%
_{l}^{(L)},\mathbf{x}_{l}^{(N)})$ (see (\ref{eq:f_z_L_cond_x_L-1})). Then, if
the observation model (\ref{eq:meas_overall_N}) and the \emph{state model}
\begin{equation}%
\begin{array}
[c]{c}%
f\left(  \mathbf{x}_{l+1}^{(N)}\left\vert \mathbf{x}_{l}^{(N)},\mathbf{x}%
_{l}^{(L)}\right.  \right) \\
=f_{\mathbf{w}_{N}}\left(  \mathbf{x}_{l+1}^{(N)}-\mathbf{f}_{l}^{(N)}\left(
\mathbf{x}_{l}^{(N)}\right)  -\mathbf{A}_{l}^{(N)}\left(  \mathbf{x}_{l}%
^{(N)}\right)  \mathbf{x}_{l}^{(L)}\right)
\end{array}
\label{eq:state_model_N-1}%
\end{equation}
are adopted for $\mathbf{x}_{l}^{(N)}$, the red graph of Fig. \ref{Fig_2} can
be drawn and exploited in a similar way as the blue graph for the evaluation
of $f(\mathbf{x}_{t}^{(N)},\mathbf{y}_{1:t},\mathbf{z}_{1:t}^{(N)})$, under
the assumption that the couple $(\mathbf{x}_{l}^{(L)},$ $\mathbf{x}%
_{l+1}^{(L)})$ is known for any $l$. Note also that, similarly to what has
been mentioned earlier about the conditional pdf $f(\mathbf{z}_{l}%
^{(L)}|\mathbf{x}_{l}^{(L)},\mathbf{x}_{l}^{(N)})$, the evaluation of
$f(\mathbf{z}_{l}^{(N)}|\mathbf{x}_{l}^{(N)})$ requires the knowledge of the
joint pdf of $\mathbf{x}_{l+1}^{(L)}$ and $\mathbf{x}_{l}^{(L)}$ conditioned
on $\mathbf{x}_{l}^{(N)}$.

Finally, merging the blue graph with the red one (i.e., adding four equality
constraint nodes for the variables $\mathbf{x}_{l}^{(L)}$ , $\mathbf{x}%
_{l}^{(N)}$, $\mathbf{x}_{l+1}^{(L)}$ and $\mathbf{x}_{l+1}^{(N)}$ shared by
the red graph and the blue one) produces the overall FG illustrated in Fig.
\ref{Fig_2}. Given this FG, we would like to follow the same line of reasoning
as that adopted for the FG of Fig. \ref{Fig_1}. In other words, given the
input messages $\vec{m}_{in}(\mathbf{x}_{l}^{(L)})=f(\mathbf{x}_{l}%
^{(L)},\mathbf{y}_{1:(l-1)})$ and $\vec{m}_{in}(\mathbf{x}_{l}^{(N)}%
)=f(\mathbf{x}_{l}^{(N)},\mathbf{y}_{1:(l-1)})$ (entering the FG along the
half edges associated with $\mathbf{x}_{l}^{(L)}$ and $\mathbf{x}_{l}^{(N)}$,
respectively), we would like to derive a \emph{forward only} message passing
algorithm based on this FG and generating the output messages $\vec{m}%
_{out}(\mathbf{x}_{l+1}^{(L)})=f(\mathbf{x}_{l+1}^{(L)},\mathbf{y}_{1:l})$ and
$\vec{m}_{out}(\mathbf{x}_{l+1}^{(N)})=f(\mathbf{x}_{l+1}^{(N)},\mathbf{y}%
_{1:l})$ (emerging from the FG along the half edges associated with
$\mathbf{x}_{l+1}^{(L)}$ and $\mathbf{x}_{l+1}^{(N)}$, respectively) on the
basis of the available a priori information and the noisy measurement
$\mathbf{y}_{l}$. Unluckily, the new FG, unlike that shown in Fig.
\ref{Fig_1}, is not cycle-free, so that any application of the SPA to it
unavoidably leads to\emph{\ approximate results} \cite{Kschischang_2001},
whatever \emph{message scheduling procedure}
\cite{Kschischang_2001,Kschischang_1998} is adopted. This consideration must
be carefully kept into account in both the derivation of MPF as a message
passing algorithm and in the development of possible modifications and
generalizations of this technique, as it will become clearer in Sections
\ref{sec:Message-Passing}-\ref{sec:Modifications-and-Extensions}.%

\section{Message Passing in Marginalized Particle
Filtering\label{sec:Message-Passing}}

In the following Section we show how the equations describing the $l$-th
recursion of MPF result from the application of the SPA to the FG shown in
Fig. \ref{Fig_2}. However, before illustrating the detailed derivation of such
equations, it is important to discuss the following relevant issues. First of
all, the MPF technique has been developed for the specific class of
GLG\emph{\ }SSMs \cite{Schon_2005,Lindsten_2016}, to which we always refer in
the following discussion. In particular, in the following we assume that: a)
$\{\mathbf{w}_{k}^{(L)}\}$ ($\{\mathbf{w}_{k}^{(N)}\}$) is a Gaussian random
process and all its elements have zero mean and covariance $\mathbf{C}%
_{w}^{(L)}$ ($\mathbf{C}_{w}^{(N)}$) for any $l$; b) $\{\mathbf{e}_{k}%
^{(L)}\}$ is a Gaussian random process having zero mean and covariance matrix
$\mathbf{C}_{e}$ for any $l$; c) all the above mentioned Gaussian processes
are statistically independent. Under these assumptions, the pdfs
$f(\mathbf{y}_{l}|\mathbf{x}_{l}^{(L)},\mathbf{x}_{l}^{(N})$, $f(\mathbf{z}%
_{l}^{(L)}|\mathbf{x}_{l}^{(L)})$ and $f(\mathbf{x}_{l+1}^{(L)}|\mathbf{x}%
_{l}^{(L)},\mathbf{x}_{l}^{(N)})$ (see (\ref{eq:f_y_L_cond_x_L-1}%
)-(\ref{eq:state_model_L-1})) are Gaussian with mean (covariance matrix)
$\mathbf{B}_{l}(\mathbf{x}_{l}^{(N)})\mathbf{x}_{l}^{(L)}+\mathbf{h}%
_{l}(\mathbf{x}_{l}^{(N)})$, $\mathbf{A}_{l}^{(N)}(\mathbf{x}_{l}%
^{(N)})\,\mathbf{x}_{l}^{(L)}$and $\mathbf{f}_{l}^{(L)}(\mathbf{x}_{l}%
^{(N)})+\mathbf{A}_{l}^{(L)}(\mathbf{x}_{l}^{(N)})\,\mathbf{x}_{l}^{(L)}$,
respectively ($\mathbf{C}_{e}$, $\mathbf{C}_{w}^{(N)}$ and $\mathbf{C}%
_{w}^{(L)}$, respectively). Similarly, the pdfs $f(\mathbf{z}_{l}%
^{(N)}|\mathbf{x}_{l}^{(N)})$ and $f(\mathbf{x}_{l+1}^{(N)}|\mathbf{x}%
_{l}^{(N)},\mathbf{x}_{l}^{(L)})$ (see (\ref{eq:meas_overall_N}) and
(\ref{eq:state_model_N-1})) are Gaussian with mean (covariance matrix)
$\mathbf{f}_{l}^{(L)}(\mathbf{x}_{l}^{(N)})$ and $\mathbf{f}_{l}%
^{(N)}(\mathbf{x}_{l}^{(N)})+\mathbf{A}_{l}^{(N)}(\mathbf{x}_{l}%
^{(N)})\,\mathbf{x}_{l}^{(L)}$, respectively ($\mathbf{C}_{w}^{(L)}$ and
$\mathbf{C}_{w}^{(N)}$, respectively).

\begin{figure}
	\centering
		\includegraphics[width=0.75\textwidth]{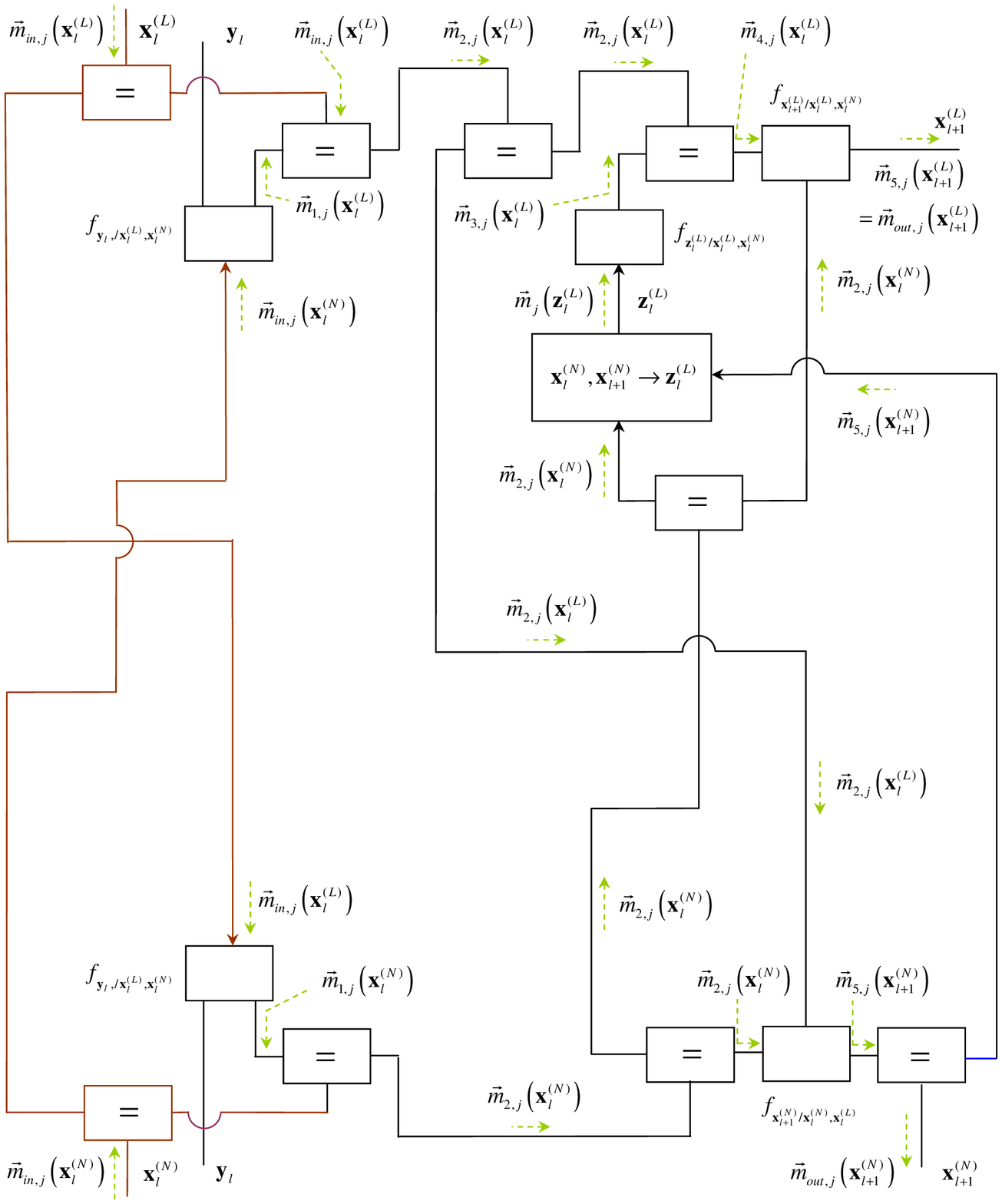}
	\caption{Overall factor graph for the representation of MPF processing; this
graph is obtained from the one shown in Fig. 3 removing the part referring to
the evaluation of \ $\mathbf{z}_{l}^{(N)}$ and inserting two new (brown)
equality constraints and some (brown) edges referring to \ $\mathbf{x}%
_{l}^{(L)}$ and $\mathbf{x}_{l}^{(N)} $. The message flow characterizing MPF
and referring to the $j$-th particle is also shown.}
	\label{Fig_3}
\end{figure}

Secondly, as explained below in detail, the MPF can be interpreted as a
\emph{forward only} message passing algorithm operating over the FG shown in
Fig. \ref{Fig_2}. The scheduling procedure adopted for MPF unavoidably leads
to ignoring the evaluation of the pseudo-measurement $\mathbf{z}_{l}^{(N)}$
(\ref{eq:z_N_l}). For this reason, in the following we refer to the simplified
FG shown in Fig. \ref{Fig_3}, which has been obtained from that illustrated in
Fig. \ref{Fig_2} removing the block representing the transformation from
$(\mathbf{x}_{l}^{(L)},\mathbf{x}_{l+1}^{(L)})$ to $\mathbf{z}_{l}^{(N)}$ and
the edges referring to the evaluation of the last vector. Note that: a) in the
new graph the block referring to the pdf $f(\mathbf{y}_{l}|\mathbf{x}%
_{l}^{(L)},\mathbf{x}_{l}^{(N})$ appears twice, since this pdf is involved in
the two subgraphs shown in Fig. \ref{Fig_2}; b) some brown edges and equality
nodes have been added to feed such blocks with $m_{in}(\mathbf{x}_{l}^{(N)})$
and $m_{in}(\mathbf{x}_{l}^{(L)})$, since these represents the only a priori
information available about $\mathbf{x}_{l}^{(N)}$ and $\mathbf{x}_{l}^{(L)}$,
respectively, at the beginning of the $l$-th recursion.

Thirdly, in MPF a particle-based model and a Gaussian model are adopted for
the input and the output messages referring to $\mathbf{x}_{l}^{(N)}$ and
$\mathbf{x}_{l}^{(L)}$, respectively, and the functional structure of the
generated messages is preserved in each recursion. More specifically, on the
one hand, the a priori information available about $\mathbf{x}_{l}^{(N)}$ at
the beginning of the $l$-th recursion is represented by a set of $N_{p}$
\emph{particles} $S_{l/(l-1)}^{(N)}=\{\mathbf{x}_{l/(l-1),j}^{(N)}%
,\,j=0,1,...,N_{p}-1\}$ and their weights $\{w_{l/(l-1),j}%
,\,j=0,\,1,...,\,N_{p}-1)\}$; following \cite{Schon_2005}, we assume that such
weights are uniform (in other words, $w_{l/(l-1),j}=1/N_{p}$ for
$j=0,\,1,...,\,N_{p}-1$), so that they can be ignored in the following
derivation. On the other hand, the a priori information available about
$\mathbf{x}_{l}^{(L)}$ is represented by a set of Gaussian pdfs, each
associated with a specific particle; in particular, the Gaussian model
$\mathcal{N}(\mathbf{x}_{l}^{(L)};\mathbf{\eta}_{l/(l.-1),j}^{(L)}%
,\mathbf{C}_{l/(l-1),j}^{(L)})$ is associated with the $j$-th particle (with
$j=0,\,1,...,\,N_{p}-1$) at the beginning of the same recursion. From the last
point it can be inferred that, in developing a message passing algorithm that
represents the MPF technique, we can focus on: a) a \emph{single particle}
contained in the input message $m_{in}(\mathbf{x}_{l}^{(N)})$ and, in
particular, on the $j$-th particle $\mathbf{x}_{l/(l-1),j}^{(N)}$; b) on the
Gaussian model $\mathcal{N}(\mathbf{x}_{l}^{(L)};\mathbf{\eta}_{l/(l.-1),j}%
^{(L)},\mathbf{C}_{l/(l-1),j}^{(L)})$ associated with that particle. For this
reason, we assume that, at the beginning of the $l$-th recursion, our
knowledge about $\mathbf{x}_{l}^{(L)}$ and $\mathbf{x}_{l}^{(N)}$ is condensed
in the message
\begin{equation}
\vec{m}_{in,j}\left(  \mathbf{x}_{l}^{(L)}\right)  =\mathcal{N}\left(
\mathbf{x}_{l}^{(L)};\mathbf{\eta}_{l/(l.-1),j}^{(L)},\mathbf{C}%
_{l/(l-1),j}^{(L)}\right) \label{eq:eq:message_L_pred-1}%
\end{equation}
and in the message
\begin{equation}
\vec{m}_{in,j}\left(  \mathbf{x}_{l}^{(N)}\right)  =\delta\left(
\mathbf{x}_{l}^{(N)}-\mathbf{x}_{l/(l-1),j}^{(N)}\right)
,\label{eq:message_N_pred-1}%
\end{equation}
respectively, with $j=0,1,...,N_{p}-1$; these are processed to generate the
corresponding output messages $\vec{m}_{out,j}(\mathbf{x}_{l+1}^{(L)})$ and
$\vec{m}_{out,j}(\mathbf{x}_{l+1}^{(N)})$, which are required to have the same
functional form as $\vec{m}_{in,j}(\mathbf{x}_{l}^{(L)})$
(\ref{eq:eq:message_L_pred-1}) and $\vec{m}_{in,j}(\mathbf{x}_{l}^{(N)})$
(\ref{eq:message_N_pred-1}), respectively. For this reason, the algorithm for
computing $\vec{m}_{out,j}(\mathbf{x}_{l+1}^{(N)})$ is expected to generate a
new particle $\mathbf{x}_{(l+1)/l,j}^{(N)}$ with a (uniform) weight
$w_{(l+1)/l,j}=1/N_{p}$; similarly, that for evaluating $\vec{m}%
_{out,j}(\mathbf{x}_{l+1}^{(L)})$ is expected to produce a new Gaussian pdf
$\mathcal{N}(\mathbf{x}_{l+1}^{(L)};\mathbf{\eta}_{(l+1)/l,j}^{(L)}%
,\mathbf{C}_{(l+1)/l,j}^{(L)})$ associated with the particle $\mathbf{x}%
_{(l+1)/l,j}^{(N)}$ (note that, in deriving this pdf, possible scale factors
are unrelevant and, consequently, can be dropped).

Given $\vec{m}_{in,j}(\mathbf{x}_{l}^{(L)})$ (\ref{eq:eq:message_L_pred-1})
and $\vec{m}_{in,j}(\mathbf{x}_{l}^{(N)})$ (\ref{eq:message_N_pred-1}), if the
SPA is applied to the considered graph and the message scheduling illustrated
in Fig. \ref{Fig_3} (and, as a matter of fact, adopted in MPF) is employed,
the steps described below are carried out to evaluate $\vec{m}_{out,j}%
(\mathbf{x}_{l+1}^{(L)})$ and $\vec{m}_{out,j}(\mathbf{x}_{l+1}^{(N)})$ in the
$l$-th recursion of MPF\footnote{In the following derivations some
mathematical results about Gassian random variables (e.g., see \cite[Par.
2.3.3]{Bishop}) and Gaussian message passing in linear models (e.g., see
\cite[Table 2, p. 1303]{Loeliger_2007}) are exploited. As far as the MPF
formulation is concerned, we always refer to that given by algorithm 1 in
\cite[Sec. II]{Schon_2005}.}.%

1. \emph{Measurement update} for $\mathbf{x}_{l}^{(N)}$ - This step aims at
updating the weight of the $j$-th particle $\mathbf{x}_{l/(l-1),j}^{(N)}$ on
the basis of the new measurements $\mathbf{y}_{l}$ (this corresponds to step
2) of algorithm 1 in \cite[Sec. II]{Schon_2005}). It involves the computation
of the messages%
\begin{equation}
\vec{m}_{1,j}\left(  \mathbf{x}_{l}^{(N)}\right)  =\int f\left(
\mathbf{y}_{l}\left\vert \mathbf{x}_{l}^{(N)},\,\mathbf{x}_{l}^{(L)}\right.
\right)  \vec{m}_{in,j}\left(  \mathbf{x}_{l}^{(L)}\right)  d\mathbf{x}%
_{l}^{(L)}\label{eq:mess_1_N_l}%
\end{equation}
and%

\begin{equation}
\vec{m}_{2,j}\left(  \mathbf{x}_{l}^{(N)}\right)  =\vec{m}_{in,j}\left(
\mathbf{x}_{l}^{(N)}\right)  \,\vec{m}_{1,j}\left(  \mathbf{x}_{l}%
^{(N)}\right)  ,\label{eq:message_2_N_l}%
\end{equation}
which provides the new importance weight for the considered particle (see Fig.
\ref{Fig_3}). Substituting the expression of $f(\mathbf{y}_{l}\mathbf{x}%
_{l}^{(N)},\,\mathbf{x}_{l}^{(L)})$ (see (\ref{eq:f_y_L_cond_x_L-1})) and
(\ref{eq:eq:message_L_pred-1}) in (\ref{eq:mess_1_N_l}) produces, after some
manipulation (see the Appendix)
\begin{equation}
\vec{m}_{1,j}\left(  \mathbf{x}_{l}^{(N)}\right)  =\mathcal{N}\left(
\mathbf{y}_{l};\mathbf{\eta}_{1,l,j}^{(N)}\left(  \mathbf{x}_{l}^{(N)}\right)
,\mathbf{C}_{1,l,j}^{(N)}\left(  \mathbf{x}_{l}^{(N)}\right)  \right)
,\label{eq:mess_y_tilde_N_l-1}%
\end{equation}
where
\begin{equation}
\mathbf{\eta}_{1,l,j}^{(N)}\left(  \mathbf{x}_{l}^{(N)}\right)  \triangleq
\mathbf{B}_{l}\left(  \mathbf{x}_{l}^{(N)}\right)  \mathbf{\eta}%
_{l/(l.-1),j}^{(L)}+\mathbf{h}_{l}\left(  \mathbf{x}_{l}^{(N)}\right)
\label{eq:eta_1b}%
\end{equation}
and%
\begin{equation}
\mathbf{C}_{1,l,j}^{(N)}\left(  \mathbf{x}_{l}^{(N)}\right)  \triangleq
\mathbf{B}_{l}\left(  \mathbf{x}_{l}^{(N)}\right)  \mathbf{C}_{l/(l-1),j}%
^{(L)}\mathbf{B}_{l}\left(  \mathbf{x}_{l}^{(N)}\right)  ^{T}+\mathbf{C}%
_{e}.\label{eq:w_1b}%
\end{equation}
Then, substituting $\vec{m}_{in,j}(\mathbf{x}_{l}^{(N)})$
(\ref{eq:message_N_pred-1}) and $\vec{m}_{1,j}(\mathbf{x}_{l}^{(N)})$
(\ref{eq:mess_y_tilde_N_l-1}) in (\ref{eq:message_2_N_l}) yields
\begin{equation}
\vec{m}_{2,j}\left(  \mathbf{x}_{l}^{(N)}\right)  =w_{l,j}\,\delta\left(
\mathbf{x}_{l}^{(N)}-\mathbf{x}_{l/(l-1),j}^{(N)}\right)
,\label{eq:message_2_N_l-1}%
\end{equation}
where\footnote{In evaluating this weight, the factor $[\det(\mathbf{C}%
_{1,l,j}^{(N)})]^{-P/2}$ appearing in the expression of the involved Gaussian
pdf \ is usually neglected, since this entails a negligible loss in estimation
accuracy.}
\begin{equation}
w_{l,j}\triangleq\mathcal{N}\left(  \mathbf{y}_{l};\mathbf{\eta}_{1,l,j}%
^{(N)},\mathbf{C}_{1,l,j}^{(N)}\right) \label{eq:weight_before_resampling}%
\end{equation}
is the new particle weight combining the a priori information about
$\mathbf{x}_{l}^{(N)}$ with the information provided by the new measurements;
here (see (\ref{eq:eta_1b}) and (\ref{eq:w_1b}))%
\begin{equation}
\mathbf{\eta}_{1,l,j}^{(N)}\triangleq\mathbf{\eta}_{1,l,j}^{(N)}\left(
\mathbf{x}_{l/(l-1),j}^{(N)}\right)  =\mathbf{B}_{l,j}\mathbf{\eta
}_{l/(l.-1),j}^{(L)}+\mathbf{h}_{l,j}\label{eq:eta_1}%
\end{equation}
and%
\begin{equation}
\mathbf{C}_{1,l,j}^{(N)}\triangleq\mathbf{C}_{1,l,j}^{(N)}\left(
\mathbf{x}_{l/(l-1),j}^{(N)}\right)  =\mathbf{B}_{l,j}\mathbf{C}%
_{l/(l-1),j}^{(L)}\mathbf{B}_{l,j}^{T}+\mathbf{C}_{e},\label{eq:w_1}%
\end{equation}
with $\mathbf{h}_{l,j}\triangleq\mathbf{h}_{l}(\mathbf{x}_{l/(l-1),j}^{(N)}) $
and $\mathbf{B}_{l,j}\triangleq\mathbf{B}_{l}(\mathbf{x}_{l/(l-1),j}^{(N)})$.
In MPF, after normalization\footnote{Note that normalization requires the
knowledge of all the weights $\{w_{l,j}\}$ (\ref{eq:weight_before_resampling})
and that, unlike it, all the previous and following tasks can be carried out
in parallel (i.e., on a particle-by-particle basis).} of the particle weights
$\{w_{l,j}\}$ (i.e., after dividing them by $P_{l}^{(w)}\triangleq%
{\displaystyle\sum\limits_{j=0}^{N_{p}-1}}
w_{l,j}$), particle \emph{resampling} with replacement is accomplished (this
corresponds to step 3) of algorithm 1 in \cite[Sec. II]{Schon_2005}). Note
that, even if resampling does not emerge from the application of SPA to the
considered graph, its use, as it will become clearer at the end of this
Section, plays an important role in the generation of the new particles for
$\mathbf{x}_{l+1}^{(N)}$. Moreover, it can be easily incorporated in our
message passing; in fact, resampling simply entails that $N_{p}$ particles
$\{\mathbf{x}_{l/(l-1),j}^{(N)}\}$ and their associated weights $\{w_{l,j}\}$
(\ref{eq:weight_before_resampling}) are replaced by the new particles
$\{\mathbf{x}_{l/l,j}^{(N)}\}$ (forming the\ new set $S_{l/l}^{(N)}$) and
their weights $\{w_{l/l,j}=1/N_{p}\}$, respectively. Consequently, $\vec
{m}_{2}(\mathbf{x}_{l}^{(N)})$ (\ref{eq:message_2_N_l-1}) is replaced by%
\begin{equation}
\vec{m}_{2,j}\left(  \mathbf{x}_{l}^{(N)}\right)  =\delta\left(
\mathbf{x}_{l}^{(N)}-\mathbf{x}_{l/l,j}^{(N)}\right)
,\label{eq:message_2_N_k-2}%
\end{equation}
since the particle weight does not depend on the index $j$. It is also worth
mentioning that, after resampling, the set of Gaussian messages $\{\vec
{m}_{in,j}(\mathbf{x}_{l}^{(L)})\}$ (\ref{eq:eq:message_L_pred-1}) needs to be
properly reordered and that the messages associated with all the discarded
particles are not propagated to the next steps.

2. \emph{Measurement update} for $\mathbf{x}_{l}^{(L)}$ - This step aims at
updating our statistical knowlege about $\mathbf{x}_{l}^{(L)}$ on the basis of
the new measurement $\mathbf{y}_{l}$ (and corresponds to step 3-a) of
algorithm 1 in \cite[Sec. II]{Schon_2005}). It involves the computation of the
messages%
\begin{equation}
\vec{m}_{1,j}\left(  \mathbf{x}_{l}^{(L)}\right)  =\int f\left(
\mathbf{y}_{l}\left\vert \mathbf{x}_{l}^{(L)},\mathbf{x}_{l}^{(N)}\right.
\right)  \vec{m}_{in,j}\left(  \mathbf{x}_{l}^{(N)}\right)  d\mathbf{x}%
_{l}^{(N)}\label{eq:message_1}%
\end{equation}
and
\begin{equation}
\vec{m}_{2,j}\left(  \mathbf{x}_{l}^{(L)}\right)  =\vec{m}_{in,j}\left(
\mathbf{x}_{l}^{(L)}\right)  \,\vec{m}_{1,j}\left(  \mathbf{x}_{l}%
^{(L)}\right)  ,\label{eq:message_4}%
\end{equation}
which represents the output of the measurement update for $\mathbf{x}%
_{l}^{(L)}$ (see Fig. \ref{Fig_3}). Substituting (\ref{eq:f_y_L_cond_x_L-1})
and (\ref{eq:message_N_pred-1}) in (\ref{eq:message_1}) produces (see
\cite[Par. 2.3.3, eq. (2.115)]{Bishop})%
\begin{equation}
\vec{m}_{1,j}\left(  \mathbf{x}_{l}^{(L)}\right)  =\mathcal{N}\left(
\mathbf{y}_{l};\mathbf{B}_{l,j}\mathbf{x}_{l}^{(L)}+\mathbf{h}_{l,j}%
,\mathbf{C}_{e}\right) \label{eq:message_1-1}%
\end{equation}
which, after some manipulation (in which unrelevant scale factors are
dropped), can be put in the Gaussian form%
\begin{equation}
\vec{m}_{1,j}\left(  \mathbf{x}_{l}^{(L)}\right)  =\mathcal{N}\left(
\mathbf{x}_{l}^{(L)};\mathbf{\eta}_{1,l,j}^{(L)},\mathbf{C}_{1,l,j}%
^{(L)}\right)  ,\label{eq:message_1_L_j}%
\end{equation}
with%
\begin{equation}
\mathbf{w}_{1,l,j}^{(L)}\triangleq\mathbf{W}_{1,l,j}^{(L)}\mathbf{\eta
}_{1,l,j}^{(L)}=\mathbf{B}_{l,j}^{T}\mathbf{W}_{e}\left(  \mathbf{y}%
_{l}-\mathbf{h}_{l,j}\right)  ,\label{eq:av_1_L_l_k-1}%
\end{equation}%
\begin{equation}
\mathbf{W}_{1,l,j}^{(L)}\triangleq\left(  \mathbf{C}_{1,l,j}^{(L)}\right)
^{-1}=\mathbf{B}_{l,j}^{T}\mathbf{W}_{e}\mathbf{B}_{l,j}%
\label{eq:cov_1_L_l_k-1}%
\end{equation}
and $\mathbf{W}_{e}\triangleq\mathbf{C}_{e}^{-1}$. Then, substituting
(\ref{eq:eq:message_L_pred-1}) and (\ref{eq:message_1_L_j}) in
(\ref{eq:message_4}) yields%
\begin{equation}%
\begin{array}
[c]{c}%
\vec{m}_{2,j}\left(  \mathbf{x}_{l}^{(L)}\right)  =\mathcal{N}\left(
\mathbf{x}_{l}^{(L)};\mathbf{\eta}_{l/(l-1),j}^{(L)},\mathbf{C}_{l/(l-1),j}%
^{(L)}\right) \\
\cdot\mathcal{N}\left(  \mathbf{x}_{l}^{(L)};\mathbf{\eta}_{1,l,j}%
^{(L)},\mathbf{C}_{1,l,j}^{(L)}\right)  ,
\end{array}
\label{eq:message_2-1}%
\end{equation}
which can be reformulated as
\begin{equation}
\vec{m}_{2,j}\left(  \mathbf{x}_{l}^{(L)}\right)  =\mathcal{N}\left(
\mathbf{x}_{l}^{(L)};\mathbf{\eta}_{2,l,j}^{(L)},\mathbf{C}_{2,l,j}%
^{(L)}\right)  ,\label{eq:message_2-l_L}%
\end{equation}
if scale factors are ignored; here,%
\begin{equation}
\mathbf{w}_{2,l,j}^{(L)}\triangleq\mathbf{W}_{2,l,j}^{(L)}\mathbf{\eta
}_{2,l,j}^{(L)}=\mathbf{w}_{l/(l-1),j}^{(L)}+\mathbf{w}_{1,l,j}^{(L)}%
,\label{eq:av_2_L_l_k}%
\end{equation}%
\begin{equation}
\mathbf{W}_{2,l,j}^{(L)}\triangleq\left(  \mathbf{C}_{2,l,j}^{(L)}\right)
^{-1}=\mathbf{W}_{l/(l-1),j}^{(L)}+\mathbf{W}_{1,l,j}^{(L)}%
,\label{eq:cov_2_L_l_k}%
\end{equation}
$\mathbf{W}_{l/(l-1),j}^{(L)}\triangleq(\mathbf{C}_{l/(l-1),j}^{(L)})^{-1}$
and $\mathbf{w}_{l/(l-1),j}^{(L)}\triangleq\mathbf{W}_{l/(l-1),j}%
^{(L)}\mathbf{\eta}_{l/(l-1),j}^{(L)}$.

3. \emph{Time update} for $\mathbf{x}_{l}^{(N)}$ - This step aims at
generating the $j$-th particle for $\mathbf{x}_{l+1}^{(N)}$ and its associated
weight (this corresponds to step 3-b) of algorithm 1 in \cite[Sec.
II]{Schon_2005}); these information are conveyed by the message $\vec{m}%
_{5,j}(\mathbf{x}_{l+1}^{(N)})$, which can be expressed as (see Fig.
\ref{Fig_3})%
\begin{equation}%
\begin{array}
[c]{c}%
\vec{m}_{5,j}\left(  \mathbf{x}_{l+1}^{(N)}\right)  =\int\int\,f\left(
\mathbf{x}_{l+1}^{(N)}\left\vert \mathbf{x}_{l}^{(L)},\mathbf{x}_{l}%
^{(N)}\right.  \right)  \,\\
\cdot\vec{m}_{2,j}\left(  \mathbf{x}_{l}^{(N)}\right)  \,\vec{m}_{2,j}\left(
\mathbf{x}_{l}^{(L)}\right)  d\mathbf{x}_{l}^{(L)}d\mathbf{x}_{l}^{(N)}.
\end{array}
\label{eq:double_integral_2}%
\end{equation}
The double integral appearing in the RHS of the last equation can be evaluated
as follows. First of all, substituting $\vec{m}_{2,j}(\mathbf{x}_{l}^{(N)})$
(\ref{eq:message_2_N_k-2}) in (\ref{eq:double_integral_2}) yields
\begin{equation}
\vec{m}_{5,j}\left(  \mathbf{x}_{l+1}^{(N)}\right)  =\int f\left(
\mathbf{x}_{l+1}^{(N)}\left\vert \mathbf{x}_{l}^{(L)},\mathbf{x}_{l/l,j}%
^{(N)}\right.  \right)  \vec{m}_{2,j}\left(  \mathbf{x}_{l}^{(L)}\right)
d\mathbf{x}_{l}^{(L)}.\label{eq:Mess_3_N}%
\end{equation}
Then, substituting the expression of $f(\mathbf{x}_{l+1}^{(N)}|\mathbf{x}%
_{l}^{(N)},\mathbf{x}_{l}^{(L)})$ (see (\ref{eq:state_model_N-1})) and
$\vec{m}_{2,j}(\mathbf{x}_{l}^{(L)})$ (\ref{eq:message_2-l_L}) in
(\ref{eq:Mess_3_N}) yields, after some manipulation, the Gaussian message (see
\cite[Par. 2.3.3, eq. (2.115)]{Bishop})
\begin{equation}
\vec{m}_{5,j}\left(  \mathbf{x}_{l+1}^{(N)}\right)  =\mathcal{N}\left(
\mathbf{x}_{l+1}^{(N)};\mathbf{\eta}_{5,l,j}^{(N)},\mathbf{C}_{5,l,j}%
^{(N)}\right)  ,\label{eq:message_5_N}%
\end{equation}
where%
\begin{equation}
\mathbf{\eta}_{5,l,j}^{(N)}\triangleq\mathbf{A}_{l,j}^{(N)}\mathbf{\eta
}_{2,l,j}^{(L)}+\mathbf{f}_{l,j}^{(N)},\label{eq:eta_5_N}%
\end{equation}%
\begin{equation}
\mathbf{C}_{5,l,j}^{(N)}\triangleq\mathbf{C}_{w}^{(N)}+\mathbf{A}_{l,j}%
^{(N)}\mathbf{C}_{2,l,j}^{(L)}\left(  \mathbf{A}_{l,j}^{(N)}\right)
^{T},\label{eq:C_5_N}%
\end{equation}
$\mathbf{A}_{l,j}^{(N)}\triangleq\mathbf{A}_{l}^{(N)}(\mathbf{x}_{l/l,j}%
^{(N)})$ and $\mathbf{f}_{l,j}^{(N)}\triangleq\mathbf{f}_{l}^{(N)}%
(\mathbf{x}_{l/l,j}^{(N)})$. Note that, in principle,
\begin{equation}
\vec{m}_{out,j}\left(  \mathbf{x}_{l+1}^{(N)}\right)  =\vec{m}_{5,j}\left(
\mathbf{x}_{l+1}^{(N)}\right)  ,\label{eq:message_5_N_l+1-1}%
\end{equation}
as it can be easily inferred from Fig. \ref{Fig_3}. However, as already
mentioned above, in MPF the output message $\vec{m}_{out}^{(j)}(\mathbf{x}%
_{l+1}^{(N)})$ is required to have the same functional form as $\vec{m}%
_{in}^{(j)}(\mathbf{x}_{l}^{(N)})$ (\ref{eq:message_N_pred-1}). This result
can be achieved a) sampling the Gaussian function $\mathcal{N}(\mathbf{x}%
_{l+1}^{(N)};\mathbf{\eta}_{5,l,j}^{(N)},\mathbf{C}_{5,l,j}^{(N)})$ (see
(\ref{eq:message_5_N})), that is drawing a sample $\mathbf{x}_{(l+1)/l,j}%
^{(N)}$ from it and b) assigning to the sample $\mathbf{x}_{(l+1)/l,j}^{(N)}$
a probability $w_{(l+1)/l,j}$ equal to the weight $w_{l/l,j}=1/N_{p}$
(originating from resampling). It is worth pointing out that: 1) the particles
$\{\mathbf{x}_{(l+1)/l,j}^{(N)}\}$ form the new set $S_{(l+1)/l}$; 2) in
accomplishing step a) of this procedure, it may be useful to introduce
artificial noise (this can be simply done adding the same positive quantity to
the diagonal elements of the matrix $\mathbf{C}_{w}^{(N)}$ appearing in the
RHS of (\ref{eq:C_5_N})) in order to mitigate the so called \emph{degeneracy
problem} \cite{Arulampalam_2002,Li_2015}. If this approach is adopted, the
message $\vec{m}_{5,j}(\mathbf{x}_{l+1}^{(N)}) $ (\ref{eq:message_5_N}) is
replaced by
\begin{equation}
\vec{m}_{5,j}\left(  \mathbf{x}_{l+1}^{(N)}\right)  =\delta\left(
\mathbf{x}_{l+1}^{(N)}-\mathbf{x}_{(l+1)/l,j}^{(N)}\right)
,\label{eq:message_5_N_l+1}%
\end{equation}
which emerges from the graph as $\vec{m}_{out,j}(\mathbf{x}_{l+1}^{(N)})$.
This message is also used in the time update for $\mathbf{x}_{l}^{(L)}$, as
illustrated in the next step.

4. \emph{Time update} for $\mathbf{x}_{l}^{(L)}$ - This step aims at
generating a new Gaussian pdf associated with the $j$-th particle
$\mathbf{x}_{(l+1)/l,j}^{(N)}$ and conveyed by $\vec{m}_{5,j}(\mathbf{x}%
_{l+1}^{(L)})=\vec{m}_{out,j}(\mathbf{x}_{l+1}^{(L)})$ (this corresponds to
step 3-c) of algorithm 1 in \cite[Sec. II]{Schon_2005}). However, before doing
that, a further measurement update is accomplished on the basis of the
pseudo-measurement $\mathbf{z}_{l}^{(L)}$ (\ref{eq:z_L_l}). This involves the
evaluation of the messages $\vec{m}_{j}(\mathbf{z}_{l}^{(L)})$,
\begin{equation}
\vec{m}_{3,j}\left(  \mathbf{x}_{l}^{(L)}\right)  =\int\vec{m}_{j}\left(
\mathbf{z}_{l}^{(L)}\right)  \,f\left(  \mathbf{z}_{l}^{(L)}\left\vert
\mathbf{x}_{l}^{(L)},\mathbf{x}_{l}^{(N)}\right.  \right)  d\mathbf{z}%
_{l}^{(L)}\label{eq:message_3}%
\end{equation}
and
\begin{equation}
\vec{m}_{4,j}\left(  \mathbf{x}_{l}^{(L)}\right)  =\vec{m}_{2,j}\left(
\mathbf{x}_{l}^{(L)}\right)  \,\vec{m}_{3,j}\left(  \mathbf{x}_{l}%
^{(L)}\right)  ,\label{eq:message_4-1}%
\end{equation}
as shown in Fig. \ref{Fig_3}. Generally speaking, the message $\vec{m}%
_{j}(\mathbf{z}_{l}^{(L)})$ can be expressed as
\begin{equation}%
\begin{array}
[c]{c}%
\vec{m}_{j}\left(  \mathbf{z}_{l}^{(L)}\right)  =\int\int f\left(
\mathbf{z}_{l}^{(L)}\left\vert \mathbf{x}_{l}^{(N)},\mathbf{x}_{l+1}%
^{(N)}\right.  \right) \\
\cdot f\left(  \mathbf{x}_{l+1}^{(N)}\left\vert \mathbf{x}_{l}^{(N)}\right.
\right)  \,f\left(  \mathbf{x}_{l}^{(N)}\right)  d\mathbf{x}_{l}%
^{(N)}d\mathbf{x}_{l+1}^{(N)}.
\end{array}
\label{eq:message_z_l}%
\end{equation}
However, since in this case $f(\mathbf{x}_{l}^{(N)})=\delta(\mathbf{x}%
_{l}^{(N)}-\mathbf{x}_{l/l,j}^{(N)})$, $f(\mathbf{x}_{l+1}^{(N)}%
|\mathbf{x}_{l}^{(L)})=\delta(\mathbf{x}_{l}^{(N+1)}-\mathbf{x}_{(l+1)/l,j}%
^{(N)})$ can be assumed (see $\vec{m}_{2,j}(\mathbf{x}_{l}^{(N)})$
(\ref{eq:message_2_N_k-2}) and $\vec{m}_{5,j}(\mathbf{x}_{l+1}^{(N)})$
(\ref{eq:message_5_N_l+1}), respectively), eq. (\ref{eq:message_z_l}) easily
leads to the expression
\begin{equation}
\vec{m}_{j}\left(  \mathbf{z}_{l}^{(L)}\right)  =f\left(  \mathbf{z}_{l}%
^{(L)}\left\vert \mathbf{x}_{l/l,j}^{(N)},\mathbf{x}_{(l+1)/l,j}^{(N)}\right.
\right)  =\delta\left(  \mathbf{z}_{l}^{(L)}-\mathbf{z}_{l,j}^{(L)}\right)
,\label{eq:message_z_L}%
\end{equation}
where
\begin{equation}
\mathbf{z}_{l,j}^{(L)}\triangleq\mathbf{x}_{(l+1)/l,j}^{(N)}-\mathbf{f}%
_{l,j}^{(N)}.\label{eq:message_Z_L}%
\end{equation}
Then, substituting (\ref{eq:message_z_L}) and the expression of $f(\mathbf{z}%
_{l}^{(L)}|\mathbf{x}_{l}^{(L)},\mathbf{x}_{l}^{(N)})$ (see
(\ref{eq:f_z_L_cond_x_L-1})) in (\ref{eq:message_3}) yields
\begin{equation}
\vec{m}_{3,j}\left(  \mathbf{x}_{l}^{(L)}\right)  =\mathcal{\mathcal{N}%
}\left(  \mathbf{z}_{l,j}^{(L)};\mathbf{A}_{l,j}^{(N)}\mathbf{x}_{l}%
^{(L)},\mathbf{C}_{w}^{(N)}\right)  .\label{eq:message_3_L}%
\end{equation}
Finally, substituting the last expression and (\ref{eq:message_2-l_L}) in
(\ref{eq:message_4-1}) produces
\begin{equation}%
\begin{array}
[c]{c}%
\vec{m}_{4,j}\left(  \mathbf{x}_{l}^{(L)}\right)  =\mathcal{N}\left(
\mathbf{x}_{l}^{(L)};\mathbf{\eta}_{2,l,j}^{(L)},\mathbf{C}_{2,l,j}%
^{(L)}\right) \\
\cdot\mathcal{\mathcal{N}}\left(  \mathbf{z}_{l,j}^{(L)};\mathbf{A}%
_{l,j}^{(N)}\mathbf{x}_{l}^{(L)},\mathbf{C}_{w}^{(N)}\right)  ,
\end{array}
\label{eq:message_4-2}%
\end{equation}
which, after some manipulation (in which unrelevant scale factors are
dropped), can be rewritten as
\begin{equation}
\vec{m}_{4,j}\left(  \mathbf{x}_{l}^{(L)}\right)  =\mathcal{N}\left(
\mathbf{x}_{l}^{(L)};\mathbf{\eta}_{4,l,j}^{(L)},\mathbf{C}_{4,l,j}%
^{(L)}\right)  ,\label{eq:message_4-3}%
\end{equation}
where%
\begin{equation}
\mathbf{w}_{4,l,j}^{(L)}\triangleq\mathbf{W}_{4,l,j}^{(L)}\mathbf{\eta
}_{4,l,j}^{(L)}=\mathbf{w}_{2,l,j}^{(L)}+\left(  \mathbf{A}_{l,j}%
^{(N)}\right)  ^{T}\mathbf{W}_{w}^{(N)}\mathbf{z}_{l,j}^{(L)}%
,\label{eq:av_4_L_l_k}%
\end{equation}%
\begin{equation}
\mathbf{W}_{4,l,j}^{(L)}\triangleq\left(  \mathbf{C}_{4,l,j}^{(L)}\right)
^{-1}=\mathbf{W}_{2,l,j}^{(L)}+\left(  \mathbf{A}_{l,j}^{(N)}\right)
^{T}\mathbf{W}_{w}^{(N)}\mathbf{A}_{l,j}^{(N)}\label{eq:cov_4_L_l_k}%
\end{equation}
and $\mathbf{W}_{w}^{(N)}\triangleq\lbrack\mathbf{C}_{w}^{(N)}]^{-1}$.

The last part of the time update step for $\mathbf{x}_{l}^{(L)}$ requires the
evaluation of the output message
\begin{equation}%
\begin{array}
[c]{c}%
\vec{m}_{5,j}\left(  \mathbf{x}_{l+1}^{(L)}\right)  =\int\int\,f\left(
\mathbf{x}_{l+1}^{(L)}\left\vert \mathbf{x}_{l}^{(L)},\mathbf{x}_{l}%
^{(N)}\right.  \right)  \,\\
\cdot\vec{m}_{4,j}\left(  \mathbf{x}_{l}^{(L)}\right)  \,\vec{m}_{2,j}\left(
\mathbf{x}_{l}^{(N)}\right)  \,d\mathbf{x}_{l}^{(L)}d\mathbf{x}_{l}^{(N)},
\end{array}
\label{eq:message_5-1}%
\end{equation}
which, similarly as $\vec{m}_{5,j}(\mathbf{x}_{l+1}^{(N)})$
(\ref{eq:double_integral_2}), requires double integration. Substituting
$\vec{m}_{2,j}(\mathbf{x}_{l}^{(N)})$ (\ref{eq:message_2_N_k-2}) in the RHS of
the last expression yields
\begin{equation}
\vec{m}_{5,j}\left(  \mathbf{x}_{l+1}^{(L)}\right)  =\int\,f\left(
\mathbf{x}_{l+1}^{(L)}\left\vert \mathbf{x}_{l}^{(L)},\mathbf{x}_{l/l,j}%
^{(N)}\right.  \right)  \,\vec{m}_{4,j}\left(  \mathbf{x}_{l}^{(L)}\right)
\,d\mathbf{x}_{l}^{(L)}.\label{eq:message_5-1-1}%
\end{equation}
Then, substituting the expression of $f(\mathbf{x}_{l+1}^{(L)}|\mathbf{x}%
_{l}^{(L)},\mathbf{x}_{l/(l-1),j}^{(N)})$ (see (\ref{eq:state_model_L-1})) and
(\ref{eq:message_4-3}) in the last equation gives (see \cite[Par. 2.3.3, eq.
(2.115)]{Bishop})
\begin{equation}%
\begin{array}
[c]{c}%
\vec{m}_{5,j}\left(  \mathbf{x}_{l+1}^{(L)}\right)  =\mathcal{N}\left(
\mathbf{x}_{l+1}^{(L)};\mathbf{\eta}_{5,l,j}^{(L)},\mathbf{C}_{5,l,j}%
^{(L)}\right) \\
=\vec{m}_{out,j}\left(  \mathbf{x}_{l+1}^{(L)}\right)  \text{,}%
\end{array}
\label{eq:message_L_out}%
\end{equation}
where
\begin{equation}
\mathbf{\eta}_{5,l,j}^{(L)}\triangleq\mathbf{A}_{l,j}^{(L)}\mathbf{\eta
}_{4,l,j}^{(L)}+\mathbf{f}_{l,j}^{(L)}=\mathbf{\eta}_{(l+1)/l,j}%
^{(L)},\label{eq:w_g_k-1-1}%
\end{equation}%
\begin{equation}
\mathbf{C}_{5,l,j}^{(L)}\triangleq\mathbf{C}_{w}^{(L)}+\mathbf{A}_{l,j}%
^{(L)}\mathbf{C}_{4,l,j}^{(L)}\left(  \mathbf{A}_{l,j}^{(L)}\right)
^{T}=\mathbf{C}_{(l+1)/l,j}^{(L)},\label{eq:W_g_k-1-1}%
\end{equation}
$\mathbf{f}_{l,j}^{(L)}\triangleq\mathbf{f}_{l}^{(L)}(\mathbf{x}_{l/l,j}%
^{(N)})$ and $\mathbf{A}_{l,j}^{(L)}\triangleq\mathbf{A}_{l}^{(L)}%
(\mathbf{x}_{l/l,j}^{(N)})$. The evaluation of $\vec{m}_{out,j}(\mathbf{x}%
_{l+1}^{(L)})$ (\ref{eq:message_L_out}) concludes the MPF message passing
procedure, which needs to be carried out for each of the $N_{p}$ particles
available at the beginning of the $l$-th recursion. Note that this procedure
needs a proper inizialization (this corresponds to step 1) of algorithm 1 in
\cite[Sec. II]{Schon_2005}). In practice, before starting the first recursion
(corresponding to $l=1$), the set $S_{1/0}^{(N)}=\{\mathbf{x}_{1/0,j}%
^{(N)},\,j=0,1,...,N_{p}-1\}$, consisting of $N_{p}$ particles, is generated
for $\mathbf{x}_{1}^{(N)}$ sampling the pdf
\begin{equation}
f\left(  \mathbf{x}_{1}^{(N)}\right)  =\int f\left(  \mathbf{x}_{1}\right)
d\mathbf{x}_{1}^{(L)},\label{eq:first_rec_N}%
\end{equation}
and the same weight $w_{1/0}=1/N_{p}$ and Gaussian model $\mathcal{N}%
(\mathbf{x}_{1}^{(L)};\mathbf{\eta}_{1/0}^{(L)},\mathbf{C}_{1/0}^{(L)})$ for
$\mathbf{x}_{1}^{(L)}$ are assigned to each of them.

Finally, it is worth pointing out that: 1) the processing accomplished in the
measurement and time update for $\mathbf{x}_{l}^{(L)}$ can be interpreted as a
form of \emph{Kalman filtering}, in which both the real measurement
$\mathbf{y}_{l}$ and the pseudo-measurement $\mathbf{z}_{l}^{(L)} $ are
processed \cite{Schon_2005}; 2) in the $l$-th recursion estimates of
$\mathbf{x}_{l}^{(N)}$ and $\mathbf{x}_{l}^{(L)}$ can be evaluated as
$\mathbf{\hat{x}}_{l}^{(N)}=\sum_{j=0}^{N_{p}-1}w_{l,j}\mathbf{x}%
_{l/(l-1),j}^{(N)}/P_{l}^{(w)}$ (see (\ref{eq:message_2_N_l-1})) and as
$\mathbf{\hat{x}}_{l}^{(L)}=\sum_{j=0}^{N_{p}-1}\mathbf{\eta}_{4,l,j}%
^{(L)}/N_{_{p}}$ (see (\ref{eq:av_4_L_l_k})), respectively; 3) the result
expressed by eq. (\ref{eq:message_5_N}) shows that, generally speaking, the
statistical representation generated by the SPA for the state $\mathbf{x}%
_{l+1}^{(N)}$ is a \emph{Gaussian mixture} (GM), whose $N_{p}$ components have
the same weight (equal to $1/N_{p}$) because resampling is always used in step
1. The last point leads to the conclusion that, if resampling was not
accomplished in the $l$-th recursion, the weight of the $j $-th component of
this GM\ would be proportional by $w_{l,j}$ (\ref{eq:weight_before_resampling}%
); this would\ unavoidably raise the problem of sampling a GM\ with unequally
weighted components in generating the particle set $S_{l/(l+1)}^{(N)}$ and
that of properly handling the resulting pseudo-measurements $\{\mathbf{z}%
_{l,j}^{(L)}\}$. These considerations motivate the use of resampling in each
recursion, indipendently of the \emph{effective sample size}
\cite{Arulampalam_2002} characterizing the particle set $S_{l/(l-1)}^{(N)}$.

\section{Simplifying Marginalized Particle Filtering\label{sec:simplifying}}

The MPF derivation illustrated in the last two Sections unveils the real
nature of MPF and its limitations, and shows the inner structure of the
processing accomplished within each step. For these reasons, it paves the way
for the development of new filtering methods related to MPF. In this Section
we exploit our FG-based representation of Bayesian filtering to develop
reduced complexity alternatives to MPF by simplifying the message passing
derived in the previous Section. It is worth mentioning that some methods for
reducing MPF computational complexity \cite{Schon_2005_complexity} have been
already proposed in the technical literature \ \cite{Smidl_2008},
\cite{Mustiere_2006}, \cite{Lu_2007}. In particular, the method proposed in
\cite{Smidl_2008} and \cite{Mustiere_2006} is based on representing the
particle set for $\mathbf{x}_{l}^{(N)}$ as a \emph{single} particle (that
corresponds to the \emph{center of mass} of the set itself), so that a
\emph{single} Kalman filter is employed in updating the statistics of the
linear component $\mathbf{x}_{l}^{(L)}$; consequently, the statistical
knowledge about $\mathbf{x}_{l}^{(L)}$ is condensed in the single message
\begin{equation}
\vec{m}_{in}\left(  \mathbf{x}_{l}^{(L)}\right)  =\mathcal{N}\left(
\mathbf{x}_{l}^{(L)};\mathbf{\eta}_{l/(l.-1)}^{(L)},\mathbf{C}_{l/(l-1)}%
^{(L)}\right)  ,\label{eq:eq:message_L_pred-1-1-1}%
\end{equation}
instead of the $N_{p}$ messages $\{\vec{m}_{in,j}(\mathbf{x}_{l}^{(L)})\}$
(\ref{eq:eq:message_L_pred-1}). Unluckily, this simplified MPF algorithm works
well only if the posterior distribution of $\mathbf{x}_{l}^{(N)}$ is
\emph{unimodal}. Its generalization to the case in which the posterior
distribution of $\mathbf{x}_{l}^{(N)}$ is \emph{multimodal} has been
illustrated later in \cite{Lu_2007}. In the proposed technique the particles
are partitioned into $K_{l}$ groups or clusters (the parameter $K_{l}$ is
required to equal the number of \emph{modes} of the posterior density of
$\mathbf{x}_{l}^{(N)}$) and each group is represented by a single particle
that corresponds to its center of mass; this allows to reduce the overall
number of Kalman filters from $N_{p}$ to $K_{l}$. The implementation of this
technique requires, however, solving the following two specific problems: a)
identifying the number of modes of the posterior distribution of
$\mathbf{x}_{l}^{(N)}$; b) partitioning the particles into clusters according
to a grouping method in each recursion. Unluckily, practical solutions for
suche problems have not been proposed in \cite{Lu_2007}.

Our derivation of simplified algorithms has been only partly inspired by the
manuscripts cited above. In fact, first of all, it relies on the following
specific methods: a) a set of $N_{p}$ equal weight particles $\{\mathbf{x}%
_{j};\,j=0,\,1,...,\,N_{p}-1\}$ is represented through its \emph{center of
mass}
\begin{equation}
\mathbf{\bar{x}}\triangleq\frac{1}{N_{p}}\sum_{j=0}^{N_{p}-1}\mathbf{x}%
_{j}\text{,}\label{center_of_mass}%
\end{equation}
as already suggested in \cite{Mustiere_2006} and \cite{Lu_2007}, when the
computation of a message referring to $\mathbf{x}_{l}^{(L)}$ involves the
particle-based representation of $\mathbf{x}_{l}^{(N)}$; b) a set of $N_{p}$
Gaussian messages $\{\mathcal{N}(\mathbf{x};\mathbf{\eta}_{j},\mathbf{C}%
_{j});\,j=0,\,1,...,\,N_{p}-1\}$, that refer to a set of $N_{p}$ equal weight
particles, is represented as the $N_{p}-$component GM
\begin{equation}
f_{GM}\left(  \mathbf{x}\right)  \triangleq\frac{1}{N_{p}}\sum_{j=0}^{N_{p}%
-1}\mathcal{N}\left(  \mathbf{x};\mathbf{\eta}_{j},\mathbf{C}_{j}\right)
,\label{eq:message_6-1}%
\end{equation}
and this GM is approximated through its projection onto the Gaussian pdf
$f_{G}\left(  \mathbf{x}\right)  =\mathcal{N}\left(  \mathbf{x};\mathbf{\eta
}_{G},\mathbf{C}_{G}\right)  $, where $\mathbf{\eta}_{G}$ and $\mathbf{C}_{G}$
are selected as%

\begin{equation}
\mathbf{\eta}_{G}\triangleq\frac{1}{N_{p}}\sum_{j=0}^{N_{p}-1}\mathbf{\eta
}_{j}\label{eq:eta_G_from_GM}%
\end{equation}
and
\begin{equation}%
\begin{array}
[t]{c}%
\mathbf{C}_{G}=(1/N_{p})\sum_{j=0}^{N_{p}-1}\mathbf{C}_{j}\\
-\mathbf{\eta}_{G}\left(  \mathbf{\eta}_{G}\right)  ^{T}+(1/N_{p})\sum
_{j=0}^{N_{p}-1}\mathbf{\eta}_{j}\left(  \mathbf{\eta}_{j}\right)  ^{T}%
\end{array}
\label{eq:cov_G_from_GM}%
\end{equation}
respectively, so that the\emph{\ mean and covariance} of $f_{GM}\left(
\mathbf{x}\right)  $ (\ref{eq:message_6-1}) are preserved (e.g., see
\cite[Sec. IV]{Runnalls_2007}). Secondly, as far as the messages $\{\vec
{m}_{in,j}(\mathbf{x}_{l}^{(L)})\}$ (\ref{eq:eq:message_L_pred-1}) are
concerned, we do not adopt the approximations proposed in \cite{Smidl_2008},
\cite{Mustiere_2006} and \cite{Lu_2007}. In fact, we focus on the following
two cases: \textbf{case \#1} - the messages $\{\vec{m}_{in,j}(\mathbf{x}%
_{l}^{(L)})\}$ are all different but, when needed in message passing, are
condensed in the single message (\ref{eq:eq:message_L_pred-1-1-1}), where
$\mathbf{\eta}_{l/(l.-1)}^{(L)}$ and $\mathbf{C}_{l/(l-1)}^{(L)}$ are given by
(\ref{eq:eta_G_from_GM}) and (\ref{eq:cov_G_from_GM}), respectively, with
$\mathbf{\eta}_{j}=\mathbf{\eta}_{l/(l.-1),j}^{(L)}$ and $\mathbf{C}%
_{j}=\mathbf{C}_{l/(l-1),j}^{(L)}$; b) \textbf{case \#2 }- the messages
$\{\vec{m}_{in,j}(\mathbf{x}_{l}^{(L)})\}$ have different means, but their
covariance matrices $\{\mathbf{C}_{l/(l-1),j}^{(L)}\}$ are all equal (their
common value is denoted $\mathbf{\tilde{C}}_{l/(l-1)}^{(L)}$ in the
following). In both cases our simplifications aim at minimizing the overall
number of a) Cholesky decompositions for the generation of the new particle
set $S_{(l+1)/l}$ (such decompositions are required for the $N_{p}$ matrices
$\{\mathbf{C}_{5,l,j}^{(L)}\}$ (\ref{eq:C_5_N})) and b) matrix inversions;
such inversions are needed to compute: a) the $N_{p}$ matrices $\{\mathbf{W}%
_{1,l,j}^{(N)}\triangleq(\mathbf{C}_{1,l,j}^{(N)})^{-1}\}$ (required in the
evaluation of the weights $\{w_{l,j}\}$ on the basis of
(\ref{eq:weight_before_resampling})); b) the $N_{p}$ matrices $\{\mathbf{W}%
_{l/(l-1),j}^{(L)}\}$ (required to evaluate the vectors $\{\mathbf{w}%
_{2,l,j}^{(L)}\}$ (\ref{eq:av_2_L_l_k}) and the matrices $\{\mathbf{W}%
_{2,l,j}^{(L)}\}$ (\ref{eq:cov_2_L_l_k})); c) the $N_{p}$ matrices
$\{\mathbf{C}_{2,l,j}^{(L)}\}$ (employed in (\ref{eq:C_5_N})); c) the $N_{p}$
matrices $\{\mathbf{C}_{4,l,j}^{(L)}\}$ (employed in (\ref{eq:W_g_k-1-1})).

Based on the methods illustrated above, our simplified versions of MPF are
derived as follows. First of all, in the \emph{measurement update for}
$\mathbf{x}_{l}^{(N)}$, we use a single covariance matrix in the Gaussian pdf
appearing in the RHS\ of (\ref{eq:weight_before_resampling}); in other words,
the $j$-th weigth $w_{l,j}$ is computed as
\begin{equation}
w_{l,j}\triangleq\mathcal{N}\left(  \mathbf{y}_{l};\mathbf{\eta}_{1,l,j}%
^{(N)},\mathbf{C}_{1,l}^{(N)}\right)  ,\label{weights_new}%
\end{equation}
where $\mathbf{C}_{1,l}^{(N)}$ is evaluated on the basis of
(\ref{eq:cov_G_from_GM}) (see also (\ref{eq:eta_G_from_GM})), setting
$\mathbf{\eta}_{j}=\mathbf{\eta}_{1,l,j}^{(N)}$ (\ref{eq:eta_1}) and
$\mathbf{C}_{j}=\mathbf{C}_{1,l,j}^{(N)}$ (\ref{eq:w_1}) for any $j$.

Second, in the \emph{measurement update for} $\mathbf{x}_{l}^{(L)}$, the
particle set $S_{l/(l-1)}^{(N)}$ is condensed in its center of mass
$\mathbf{\bar{x}}_{l/(l-1)}^{(N)}$ (see (\ref{center_of_mass})). Consequently,
the message $\vec{m}_{1,j}(\mathbf{x}_{l}^{(L)})$ (\ref{eq:message_1-1}) is
replaced by its particle-independent form%
\begin{equation}
\vec{m}_{1}\left(  \mathbf{x}_{l}^{(L)}\right)  =\mathcal{N}\left(
\mathbf{y}_{l};\mathbf{\bar{B}}_{l}\mathbf{x}_{l}^{(L)}+\mathbf{\bar{h}}%
_{l},\mathbf{C}_{e}\right) \label{message_1_L_av}%
\end{equation}
where $\mathbf{\bar{B}}_{l}\triangleq\mathbf{B}_{l}(\mathbf{\bar{x}}%
_{l/(l-1)}^{(N)})$ and $\mathbf{\bar{h}}_{l}\triangleq\mathbf{h}%
_{l}(\mathbf{\bar{x}}_{l/(l-1)}^{(N)})$. This message, similarly as
(\ref{eq:message_1-1}), can be put in the Gaussian form (see
(\ref{eq:message_1_L_j})-(\ref{eq:cov_1_L_l_k-1}))%
\begin{equation}
\vec{m}_{1}\left(  \mathbf{x}_{l}^{(L)}\right)  =\mathcal{N}\left(
\mathbf{x}_{l}^{(L)};\mathbf{\eta}_{1,l}^{(L)},\mathbf{C}_{1,l}^{(L)}\right)
,\label{message_1_L_av_bis}%
\end{equation}
where%
\begin{equation}
\mathbf{w}_{1,l}^{(L)}\triangleq\mathbf{W}_{1,l}^{(L)}\mathbf{\eta}%
_{1,l}^{(L)}=\mathbf{\bar{B}}_{l}^{T}\mathbf{W}_{e}\left(  \mathbf{y}%
_{l}-\mathbf{\bar{h}}_{l}\right) \label{message_1_L_av_w}%
\end{equation}
and%
\begin{equation}
\mathbf{W}_{1,l}^{(L)}\triangleq\left(  \mathbf{C}_{1,l}^{(L)}\right)
^{-1}=\mathbf{\bar{B}}_{l}^{T}\mathbf{W}_{e}\mathbf{\bar{B}}%
.\label{message_1_L_av_W}%
\end{equation}
This allows us to replace the message $\vec{m}_{2,j}(\mathbf{x}_{l}^{(L)})$
(\ref{eq:message_2-l_L}) with its particle-independent counterpart
\begin{equation}
\vec{m}_{2}\left(  \mathbf{x}_{l}^{(L)}\right)  =\mathcal{N}\left(
\mathbf{x}_{l}^{(L)};\mathbf{\eta}_{2,l}^{(L)},\mathbf{C}_{2,l}^{(L)}\right)
,\label{eq:message_2-l_L_new}%
\end{equation}
where $\mathbf{w}_{2,l}^{(L)}\triangleq\mathbf{W}_{2,l}^{(L)}\mathbf{\eta
}_{2,l}^{(L)}$ and $\mathbf{W}_{2,l}^{(L)}\triangleq(\mathbf{C}_{2,l}%
^{(L)})^{-1}$ are easily obtained from (\ref{eq:av_2_L_l_k}) and
(\ref{eq:cov_2_L_l_k}) replacing a) $\mathbf{w}_{1,l,j}^{(L)}$ \ and
$\mathbf{W}_{1,l,j}^{(L)}$ with $\mathbf{w}_{1,l}^{(L)}$
(\ref{message_1_L_av_w}) and $\mathbf{W}_{1,l}^{(L)}$ (\ref{message_1_L_av_W}%
), respectively; b) $\mathbf{w}_{l/(l-1),j}^{(L)}$ and $\mathbf{W}%
_{l/(l-1),j}^{(L)}$ with $\mathbf{w}_{l/(l-1)}^{(L)}\triangleq\mathbf{W}%
_{l/(l-1)}^{(L)}\mathbf{\eta}_{l/(l-1)}^{(L)}$ and $\mathbf{W}_{l/(l-1)}%
^{(L)}\triangleq(\mathbf{C}_{l/(l-1)}^{(L)})^{-1}$ ($\mathbf{\eta}%
_{l/(l.-1)}^{(L)}$ and $\mathbf{C}_{l/(l-1)}^{(L)}$ are given by
(\ref{eq:eta_G_from_GM}) and (\ref{eq:cov_G_from_GM}), respectively, with
$\mathbf{\eta}_{j}=\mathbf{\eta}_{l/(l.-1),j}^{(L)}$ and $\mathbf{C}%
_{j}=\mathbf{C}_{l/(l-1),j}^{(L)}$). Note that, since the precision matrix
$\mathbf{W}_{2,l}^{(L)}$ is particle-independent, a \emph{single} matrix
$\mathbf{C}_{2,l}^{(L)}$ has to be computed for the next step.

Thirdly, in the \emph{time update for} $\mathbf{x}_{l}^{(N)}$, the message
$\vec{m}_{2}(\mathbf{x}_{l}^{(L)})$ (\ref{eq:message_2-l_L_new}) can be used
in place of $\vec{m}_{2,j}(\mathbf{x}_{l}^{(L)})$ (\ref{eq:message_2-l_L}) in
the evaluation of $\vec{m}_{5,j}(\mathbf{x}_{l+1}^{(N)})$ (see
(\ref{eq:Mess_3_N}) and Fig. \ref{Fig_3}). This leads to the
particle-dependent message
\begin{equation}
\vec{m}_{5,j}\left(  \mathbf{x}_{l+1}^{(N)}\right)  =\mathcal{N}\left(
\mathbf{x}_{l+1}^{(N)};\mathbf{\eta}_{5,l,j}^{(N)},\mathbf{C}_{5,l,j}%
^{(N)}\right)  ,\label{eq:message_5_N-1}%
\end{equation}
where $\mathbf{\eta}_{5,l,j}^{(N)}$ and $\mathbf{C}_{5,l,j}^{(N)}$ are
obtained from (\ref{eq:eta_5_N}) and (\ref{eq:C_5_N}), respectively, replacing
$\mathbf{\eta}_{2,l,j}^{(L)}$ and $\mathbf{C}_{2,l,j}^{(L)}$ with
$\mathbf{\eta}_{2,l}^{(L)}$ and $\mathbf{C}_{2,l}^{(L)}$, respectively. Then,
to simplify the generation of the new particle set $S_{(l+1)/l}^{(N)}$, the
covariance matrices $\{\mathbf{C}_{5,l,j}^{(N)}\}$ are condensed in a single
matrix $\mathbf{C}_{5,l}^{(N)}$ using (\ref{eq:cov_G_from_GM}) (see also
(\ref{eq:eta_G_from_GM})) with $\mathbf{\eta}_{j}=\mathbf{\eta}_{5,l,j}%
^{(N)}=\mathbf{A}_{l,j}^{(N)}\mathbf{\eta}_{2,l}^{(L)}+\mathbf{f}_{l,j}^{(N)}$
and $\mathbf{C}_{j}=\mathbf{C}_{5,l,j}^{(N)}=\mathbf{C}_{w}^{(N)}%
+\mathbf{A}_{l,j}^{(N)}\mathbf{C}_{2,l}^{(L)}(\mathbf{A}_{l,j}^{(N)})^{T}$;
consequently, the particle generation mechanism for $\mathbf{x}_{l+1}^{(N)}$
requires computing the Cholesky decomposition of a \emph{single} matrix
(namely, $\mathbf{C}_{5,l}^{(N)}$), since it is based on the modified message
\begin{equation}
\vec{m}_{5}\left(  \mathbf{x}_{l+1}^{(N)}\right)  =\mathcal{N}\left(
\mathbf{x}_{l+1}^{(N)};\mathbf{\eta}_{5,l,j}^{(N)},\mathbf{C}_{5,l}%
^{(N)}\right)  ,\label{eq:message_5_N-2}%
\end{equation}
which depends on the particle index $j$ through the mean vector $\mathbf{\eta
}_{5,l,j}^{(N)}$ only.

Finally, in the \emph{time update for} $\mathbf{x}_{l}^{(L)}$, the pdf
\begin{equation}
f\left(  \mathbf{z}_{l}^{(L)}\left\vert \mathbf{x}_{l}^{(L)}\right.  \right)
=\mathcal{\mathcal{N}}\left(  \mathbf{z}_{l}^{(L)};\mathbf{\bar{A}}_{l}%
^{(N)}\mathbf{x}_{l}^{(L)},\mathbf{C}_{w}^{(N)}\right) \label{pdf_z}%
\end{equation}
is employed in the evaluation of $\vec{m}_{3,j}(\mathbf{x}_{l}^{(L)})$ through
(\ref{eq:message_3}), where $\mathbf{\bar{A}}_{l}^{(N)}=\mathbf{A}_{l}%
^{(N)}(\mathbf{\bar{x}}_{l/l}^{(N)})$ and $\mathbf{\bar{x}}_{l/l}^{(N)}$
denotes the center of mass of the particle set $S_{l/l}^{(N)}$ (see
(\ref{center_of_mass})). Then, the messages $\vec{m}_{3,j}(\mathbf{x}%
_{l}^{(L)})$ (\ref{eq:message_3_L}) and $\vec{m}_{4,j}(\mathbf{x}_{l}^{(L)})$
(\ref{eq:message_4-3}) can be replaced by
\begin{equation}
\vec{m}_{3,j}\left(  \mathbf{x}_{l}^{(L)}\right)  =\mathcal{\mathcal{N}%
}\left(  \mathbf{z}_{l,j}^{(L)};\mathbf{\bar{A}}_{l}^{(N)}\mathbf{x}_{l}%
^{(L)},\mathbf{C}_{w}^{(N)}\right) \label{eq:message_3_L-2}%
\end{equation}
and
\begin{equation}
\vec{m}_{4,j}\left(  \mathbf{x}_{l}^{(L)}\right)  =\mathcal{N}\left(
\mathbf{x}_{l}^{(L)};\mathbf{\eta}_{4,l,j}^{(L)},\mathbf{C}_{4,l}%
^{(L)}\right)  ,\label{eq:message_4-2-new}%
\end{equation}
respectively, with (see (\ref{eq:av_4_L_l_k}) and (\ref{eq:cov_4_L_l_k}))%
\begin{equation}
\mathbf{w}_{4,l,j}^{(L)}\triangleq\mathbf{W}_{4,l}^{(L)}\mathbf{\eta}%
_{4,l,j}^{(L)}=\mathbf{w}_{2,l}^{(L)}+\left(  \mathbf{\bar{A}}_{l}%
^{(N)}\right)  ^{T}\mathbf{W}_{w}^{(N)}\mathbf{z}_{l,j}^{(L)}%
\label{eq:av_4_L_l_k-2}%
\end{equation}
and
\begin{equation}
\mathbf{W}_{4,l}^{(L)}\triangleq\left(  \mathbf{C}_{4,l}^{(L)}\right)
^{-1}=\mathbf{W}_{2,l}^{(L)}+\left(  \mathbf{\bar{A}}_{l}^{(N)}\right)
^{T}\mathbf{W}_{w}^{(N)}\mathbf{\bar{A}}_{l}^{(N)}.\label{eq:cov_4_L_l_k-2}%
\end{equation}
Substituting now $\vec{m}_{4,j}(\mathbf{x}_{l}^{(L)})$
(\ref{eq:message_4-2-new}) (in place of $\vec{m}_{4,j}(\mathbf{x}_{l}^{(L)})$
(\ref{eq:message_4-3})) and $\vec{m}_{2,j}(\mathbf{x}_{l}^{(N)})$
(\ref{eq:message_2_N_k-2}) in (\ref{eq:message_5-1}), produces, after some
manipulation
\begin{equation}
\vec{m}_{5,j}\left(  \mathbf{x}_{l+1}^{(L)}\right)  =\mathcal{N}\left(
\mathbf{x}_{l+1}^{(L)};\mathbf{\eta}_{5,l,j}^{(L)},\mathbf{C}_{5,l,j}%
^{(L)}\right)  ,\label{eq:message_5_new}%
\end{equation}
where
\begin{equation}
\mathbf{\eta}_{5,l,j}^{(L)}\triangleq\mathbf{A}_{l,j}^{(L)}\mathbf{\eta
}_{4,l,j}^{(L)}+\mathbf{f}_{l,j}^{(L)}\label{eq:eta_5_new}%
\end{equation}
and
\begin{equation}
\mathbf{C}_{5,l,j}^{(L)}\triangleq\mathbf{C}_{w}^{(L)}+\mathbf{A}_{l,j}%
^{(L)}\mathbf{C}_{4,l}^{(L)}\left(  \mathbf{A}_{l,j}^{(L)}\right)
^{T}.\label{eq:C_5_new}%
\end{equation}
Note that: a) since the precision matrix $\mathbf{W}_{4,l}^{(L)}$
(\ref{eq:cov_4_L_l_k-2}) is particle-independent, a \emph{single} matrix
inversion is needed to evaluate the matrix $\mathbf{C}_{4,l}^{(L)}$ appearing
in (\ref{eq:C_5_new}); b) in case \# 2 the matrix set $\{\mathbf{C}%
_{5,l,j}^{(L)}\}$ is condensed in a single matrix $\mathbf{C}_{5,l}^{(L)}$
(this represents the common value $\mathbf{\tilde{C}}_{(l+1)/l}^{(L)}$ taken
on by all the matrices $\{\mathbf{C}_{(l+1)/l,j}^{(L)}\}$ processed in the
next recursion); c) the matrix $\mathbf{C}_{5,l}^{(L)}$ is evaluated on the
basis of (\ref{eq:eta_G_from_GM}) and (\ref{eq:cov_G_from_GM}), setting
$\mathbf{\eta}_{j}=\mathbf{\eta}_{5,l,j}^{(L)}$ and $\mathbf{C}_{j}%
=\mathbf{C}_{5,l,j}^{(L)}$ for any $j$.

The new filtering techniques, based on MPF and on the simplified messages
derived above, are called \emph{simplified} MPF (SMPF) in the following; in
particular the acronyms SMPF1 and SMPF2 are used to refer to case \#1 and case
\#2, respectively.

\section{Message Passing in Iterative Filtering Techniques Inspired by
Marginalized Particle Filtering\label{sec:Modifications-and-Extensions}}

As already mentioned above, the suboptimality of MPF can related to the fact
that the FG underlying the considered filtering problem is not cycle free. It
is well known that the SPA can also be applied to a factor graph with cycles
simply by following the same message propagation rules; however, generally
speaking, this leads to an \textquotedblleft iterative\textquotedblright{}
algorithm with no natural termination (and known as \emph{loopy belief
propagation}), since its messages are passed multiple times on a given edge
\cite{Loeliger_2007}, \cite{Kschischang_2001}. Despite this, some of the most
relevant applications of the SPA have been developed for systems in which the
underlying FG does have cycles, like the one shown in Fig. \ref{Fig_1}. In the
following we show how a novel iterative technique can be developed for our
filtering problem following this approach. To begin, we note that our interest
in iterative methods is also motivated by the possibility of exploiting the
pseudo-measurement $\mathbf{z}_{l}^{(N)}$ (\ref{eq:z_N_l}); in fact, the
message $m_{j}(\mathbf{z}_{l}^{(N)})$ referring to this random vector cannot
be computed in MPF because of the adopted scheduling, but can certainly
provide additional information to refine our statistical knowledge about the
nonlinear component $\mathbf{x}_{l}^{(N)}$. This message, similarly as the one
referring to\footnote{If $\mathbf{x}_{l}^{(N)}=\mathbf{x}_{l/l,j}^{(N)}$, the
random vector $\mathbf{z}_{l}^{(L)}$ (\ref{eq:z_L_l}) becomes $\mathbf{x}%
_{l+1}^{(N)}-\mathbf{f}_{l}^{(N)}(\mathbf{x}_{l/l,j}^{(N)})$; consequently,
adopting the Gaussian model (\ref{eq:message_5_N}) for $\mathbf{x}_{l+1}%
^{(N)}$ results in a Gaussian model for $\mathbf{z}_{l}^{(L)}$ too.}
$\mathbf{z}_{l}^{(L)}$, can be put in a Gaussian form, that is
\begin{equation}
\vec{m}_{j}\left(  \mathbf{z}_{l}^{(N)}\right)  =\mathcal{N}\left(
\mathbf{z}_{l}^{(N)};\mathbf{\eta}_{\mathbf{z},l,j}^{(N)},\mathbf{C}%
_{\mathbf{z},l,j}^{(N)}\right)  ,\label{eq:message_z_N1}%
\end{equation}
since $\mathbf{x}_{l+1}^{(L)}$ and $\mathbf{x}_{l}^{(L)}$, conditioned on
$\mathbf{x}_{l}^{(N)}$, are modelled as jointly Gaussian random vectors. Let
us show now how this message can be computed in an iterative filtering
procedure generalising MPF and how it can exploited in such a procedure. First
of all, we assume that the message $\vec{m}_{5,j}(\mathbf{x}_{l+1}^{(L)})$
(\ref{eq:message_L_out}), representing the pdf of $\mathbf{x}_{l+1}^{(L)}$
conditioned on $\mathbf{x}_{l}^{(N)}=\mathbf{x}_{l/l,j}^{(N)}$, is already
available when the time update for $\mathbf{x}_{l}^{(L)}$ (i.e., the first
step of MPF) is accomplished. Then, given $\vec{m}_{2,j}(\mathbf{x}_{l}%
^{(L)})$ (\ref{eq:message_2-l_L}) and $\vec{m}_{5,j}(\mathbf{x}_{l+1}^{(L)})$
(\ref{eq:message_L_out}), the mean and covariance of $\mathbf{z}_{l}^{(N)}$
can be evaluated as (see (\ref{eq:z_N_l}))%

\begin{equation}
\mathbf{\eta}_{\mathbf{z},l,j}^{(N)}=\mathbf{\eta}_{5,l,j}^{(L)}%
-\mathbf{A}_{l,j}^{(L)}\mathbf{\eta}_{2,l,j}^{(L)}\label{eq:av_z}%
\end{equation}
and
\begin{equation}%
\begin{array}
[c]{c}%
\mathbf{C}_{\mathbf{z},l,j}^{(N)}=\mathbf{C}_{5,l,j}^{(L)}+\mathbf{A}%
_{l,j}^{(L)}\mathbf{C}_{2,l,j}^{(L)}\left(  \mathbf{A}_{l,j}^{(L)}\right)
^{T}\\
-\mathbf{A}_{l,j}^{(L)}\mathbf{C}_{\mathbf{x},l,j}^{(L)}-\left(
\mathbf{C}_{\mathbf{x},l,j}^{(L)}\right)  ^{T}\left(  \mathbf{A}_{l,j}%
^{(L)}\right)  ^{T},
\end{array}
\label{eq:C_z}%
\end{equation}
respectively, where $\mathbf{C}_{\mathbf{x},l,j}^{(L)}$ denotes the cross
covariance matrix for the vectors $\mathbf{x}_{l}^{(L)}$ and $\mathbf{x}%
_{l+1}^{(L)}$ (conditioned on $\mathbf{x}_{l}^{(N)}=\mathbf{x}_{l/l,j}^{(N)}%
$). Given $\vec{m}_{2,j}(\mathbf{x}_{l}^{(L)})$ (\ref{eq:message_2-l_L}) and
the conditional pdf $f(\mathbf{x}_{l+1}^{(L)}|\mathbf{x}_{l}^{(L)}%
,\mathbf{x}_{l/l,j}^{(N)})=\mathcal{N}(\mathbf{x}_{l+1}^{(L)};\mathbf{f}%
_{l,j}^{(L)}+\mathbf{A}_{l,j}^{(L)}\mathbf{x}_{l}^{(L)},\mathbf{C}_{w}^{(L)})$
(see (\ref{eq:state_model_L-1})), it is easy to show that $\mathbf{C}%
_{\mathbf{x},l,j}^{(L)}=\mathbf{C}_{2,l,j}^{(L)}(\mathbf{A}_{l,j}^{(L)})^{T}$
(e.g., see \cite[Par. 2.3.3, eq. (2.104)]{Bishop}); consequently, eq.
(\ref{eq:C_z}) can be rewritten as
\begin{equation}
\mathbf{C}_{\mathbf{z},l,j}^{(N)}=\mathbf{C}_{5,l,j}^{(L)}-\mathbf{A}%
_{l,j}^{(L)}\mathbf{C}_{2,l,j}^{(L)}\left(  \mathbf{A}_{l,j}^{(L)}\right)
^{T}.\label{eq:C_z-1}%
\end{equation}
Equations (\ref{eq:av_z}) and (\ref{eq:C_z-1}) represent the desired result,
since they provide a complete statistical characterization of the message
$\vec{m}_{j}(\mathbf{z}_{l}^{(N)})$ (\ref{eq:message_z_N1}). In principle,
this message could be exploited in a similar way as that adopted for $\vec
{m}_{j}(\mathbf{z}_{l}^{(L)})$ (\ref{eq:message_z_L}); this approach would
lead to draw a set of $N_{p}$ samples $\{\mathbf{z}_{l,j}^{(N)}\}$ from the
Gaussian function appearing in the RHS of (\ref{eq:message_z_N1}) and to
process the resulting pseudo-measurements to generate a new weight for each
particle of the set $S_{l/l}$. However, our computer simulations have shown
that this approach is outperformed by more refined method illustrated in the
following. In practice, in the iterative filtering method we propose the
message $\vec{m}_{j}(\mathbf{z}_{l}^{(N)})$ (\ref{eq:message_z_N1}), similarly
as $\vec{m}_{j}(\mathbf{z}_{l}^{(L)})$ (\ref{eq:message_z_L}), is employed to
evaluate the new message\footnote{Note that the following message represents
the \emph{correlation} between the pdf $\vec{m}_{j}(\mathbf{z}_{l}^{(N)})$
evaluated on the basis of the definition of $\mathbf{z}_{l}^{(N)}$
(\ref{eq:z_N_l}) and the pdf originating from the fact that this quantity is
expected to equal the random vector $\mathbf{f}_{l}^{(L)}(\mathbf{x}%
_{l/l,j}^{(N)})+\mathbf{w}_{l}^{(L)}$. For this reason, it expresses the
\emph{degree of similarity} between these two functions.}
\begin{equation}
\vec{m}_{3,j}\left(  \mathbf{x}_{l}^{(N)}\right)  =\int\vec{m}_{j}\left(
\mathbf{z}_{l}^{(N)}\right)  f\left(  \mathbf{z}_{l}^{(N)}\left\vert
\mathbf{x}_{l/l,j}^{(N)}\right.  \right)  d\mathbf{z}_{l}^{(N)}%
,\label{eq:message_3-1}%
\end{equation}
which represents for $\mathbf{x}_{l}^{(N)}$ the counterpart of the message
$\vec{m}_{3,j}(\mathbf{x}_{l}^{(L)})$ (\ref{eq:message_3}). Substituting
(\ref{eq:message_z_N1}) and the expression of $f(\mathbf{z}_{l}^{(N)}%
|\mathbf{x}_{l}^{(N)})$ (given $\mathbf{x}_{l}^{(N)}=\mathbf{x}_{l/l,j}^{(N)}%
$) in the RHS of the last expression gives\footnote{In our computer
simulations the factor $D_{3,l,j}^{(N)}$ appearing in this weight has been
always neglected, since it negligibly influences estimation accuracy.} (see
the Appendix)%
\begin{equation}%
\begin{array}
[c]{c}%
\vec{m}_{3,j}\left(  \mathbf{x}_{l}^{(N)}\right) \\
=D_{3,l,j}^{(N)}\cdot\exp\left[  \frac{1}{2}\left(  \left(  \mathbf{\eta
}_{3,l,j}^{(N)}\right)  ^{T}\mathbf{W}_{3,l,j}^{(N)}\mathbf{\eta}%
_{3,l,j}^{(N)}\right.  \right. \\
\left.  \left.  -\left(  \mathbf{\eta}_{\mathbf{z},l,j}^{(N)}\right)
^{T}\mathbf{W}_{\mathbf{z},l,j}^{(N)}\mathbf{\eta}_{\mathbf{z},l,j}%
^{(N)}-\left(  \mathbf{f}_{l,j}^{(L)}\right)  ^{T}\mathbf{W}_{w}%
^{(L)}\mathbf{f}_{l,j}^{(L)}\right)  \right]  \triangleq p_{l,j},
\end{array}
\label{eq:P_1_l_j-1}%
\end{equation}
where%
\begin{equation}
\mathbf{W}_{3,l,j}^{(N)}\triangleq\left(  \mathbf{C}_{3,l,j}^{(N)}\right)
^{-1}=\mathbf{W}_{\mathbf{z},l,j}^{(N)}+\mathbf{W}_{w}^{(L)},\label{eq:w_1-1}%
\end{equation}%
\begin{equation}
\mathbf{w}_{3,l,j}^{(N)}\triangleq\mathbf{W}_{3,l,j}^{(N)}\mathbf{\eta
}_{3,l,j}^{(N)}=\mathbf{w}_{\mathbf{z},l,j}^{(N)}+\mathbf{W}_{w}%
^{(L)}\mathbf{f}_{l,j}^{(L)},\label{eq:eta_1-1}%
\end{equation}
$\mathbf{W}_{\mathbf{z},l,j}^{(N)}\triangleq(\mathbf{C}_{\mathbf{z},l,j}%
^{(N)})^{-1}$, $\mathbf{w}_{\mathbf{z},l,j}^{(N)}\triangleq\mathbf{W}%
_{\mathbf{z},l,j}^{(N)}\mathbf{\eta}_{\mathbf{z},l,j}^{(N)}$, $\mathbf{W}%
_{w}^{(L)}\triangleq\lbrack\mathbf{C}_{w}^{(L)}]^{-1}$,%
\begin{equation}
D_{3,l,j}^{(N)}\triangleq\left[  \det\left(  \mathbf{\tilde{C}}_{l,j}%
^{(N)}\right)  \right]  ^{-D_{L}/2}\label{eq:mess_1_N_l-2-2-2-1}%
\end{equation}
and $\mathbf{\tilde{C}}_{l,j}^{(N)}\triangleq\mathbf{C}_{\mathbf{z},l,j}%
^{(N)}+\mathbf{C}_{w}^{(L)}$. Then, the new message $\vec{m}_{3,j}%
(\mathbf{x}_{l}^{(N)})$ (\ref{eq:P_1_l_j-1}) is exploited, similarly as
$\vec{m}_{3,j}(\mathbf{x}_{l}^{(L)})$ (\ref{eq:message_3_L}), to generate the
message (see (\ref{eq:message_4-1}))
\begin{equation}
\vec{m}_{4,j}\left(  \mathbf{x}_{l}^{(N)}\right)  =\vec{m}_{2,j}\left(
\mathbf{x}_{l}^{(N)}\right)  \vec{m}_{3,j}\left(  \mathbf{x}_{l}^{(N)}\right)
,\label{eq:message_4-1-1}%
\end{equation}
where $\vec{m}_{2,j}(\mathbf{x}_{l}^{(N)})$ is expressed by
(\ref{eq:message_2_N_l-1}) (i.e., it is the message emerging from the
measurement update for $\mathbf{x}_{l}^{(N)}$ in the absence of resampling).
This produces (see (\ref{eq:P_1_l_j-1}))
\begin{equation}
\vec{m}_{4,j}\left(  \mathbf{x}_{l}^{(N)}\right)  =W_{l,j}\,\delta\left(
\mathbf{x}_{l}^{(N)}-\mathbf{x}_{l/(l-1),j}^{(N)}\right)
,\label{eq:message_4_N_l-1-1}%
\end{equation}
where
\begin{equation}
W_{l,j}\triangleq w_{l,j}\cdot p_{l,j}\label{eq:weight_before_resampling-1}%
\end{equation}
represents the new weight for the $j$-th particle; such a weight accounts for
both the (real) measurement $\mathbf{y}_{l}$ and the pseudo-measurement
$\mathbf{z}_{l,j}^{(N)}$. Resampling with replacement can now be accomplished
for the set $S_{l/(l-1)}^{(N)}\triangleq\{\mathbf{x}_{l/(l-1),j}^{(N)}\}$ on
the basis of the more refined weights $\{W_{l,j}\}$
(\ref{eq:weight_before_resampling-1}); if this is done, the message $\vec
{m}_{4,j}(\mathbf{x}_{l}^{(N)})$ takes on the same form as $\vec{m}%
_{2,j}(\mathbf{x}_{l}^{(N)})$ (\ref{eq:message_2_N_k-2}), i.e. it can be
expressed as
\begin{equation}
\vec{m}_{4,j}\left(  \mathbf{x}_{l}^{(N)}\right)  =\delta\left(
\mathbf{x}_{l}^{(N)}-\mathbf{\tilde{x}}_{l/l,j}^{(N)}\right)
,\label{eq:message_4_N_l-1-1-1}%
\end{equation}
with $j=0,1,...,N_{p}-1$. Note that, generally speaking, the particle set
$\tilde{S}_{l/l}^{(N)}\triangleq\{\mathbf{\tilde{x}}_{l/l,j}^{(N)}\}$ produced
by resampling in this case is different from $S_{l/l}^{(N)}\triangleq
\{\mathbf{x}_{l/l,j}^{(N)}\}$ (i.e., from the one obtained with MPF), even if
both of them originate from the same set $S_{l/(l-1)}^{(N)}$; this is due to
the fact that the weights $\{W_{l,j}\}$ (\ref{eq:weight_before_resampling-1})
may be substantially different from the MPF weights $\{w_{l,j}\}$ because of
the factor $p_{l,j}$. Finally, the new message $\vec{m}_{4,j}(\mathbf{x}%
_{l}^{(N)})$ (\ref{eq:message_4_N_l-1-1-1}) is used in place of $\vec{m}%
_{2,j}(\mathbf{x}_{l}^{(N)})$ in the RHS of (\ref{eq:double_integral_2}) for
the evaluation of the message $\vec{m}_{5,j}(\mathbf{x}_{l+1}^{(L)})$.

\begin{figure}
	\centering
		\includegraphics[width=0.75\textwidth]{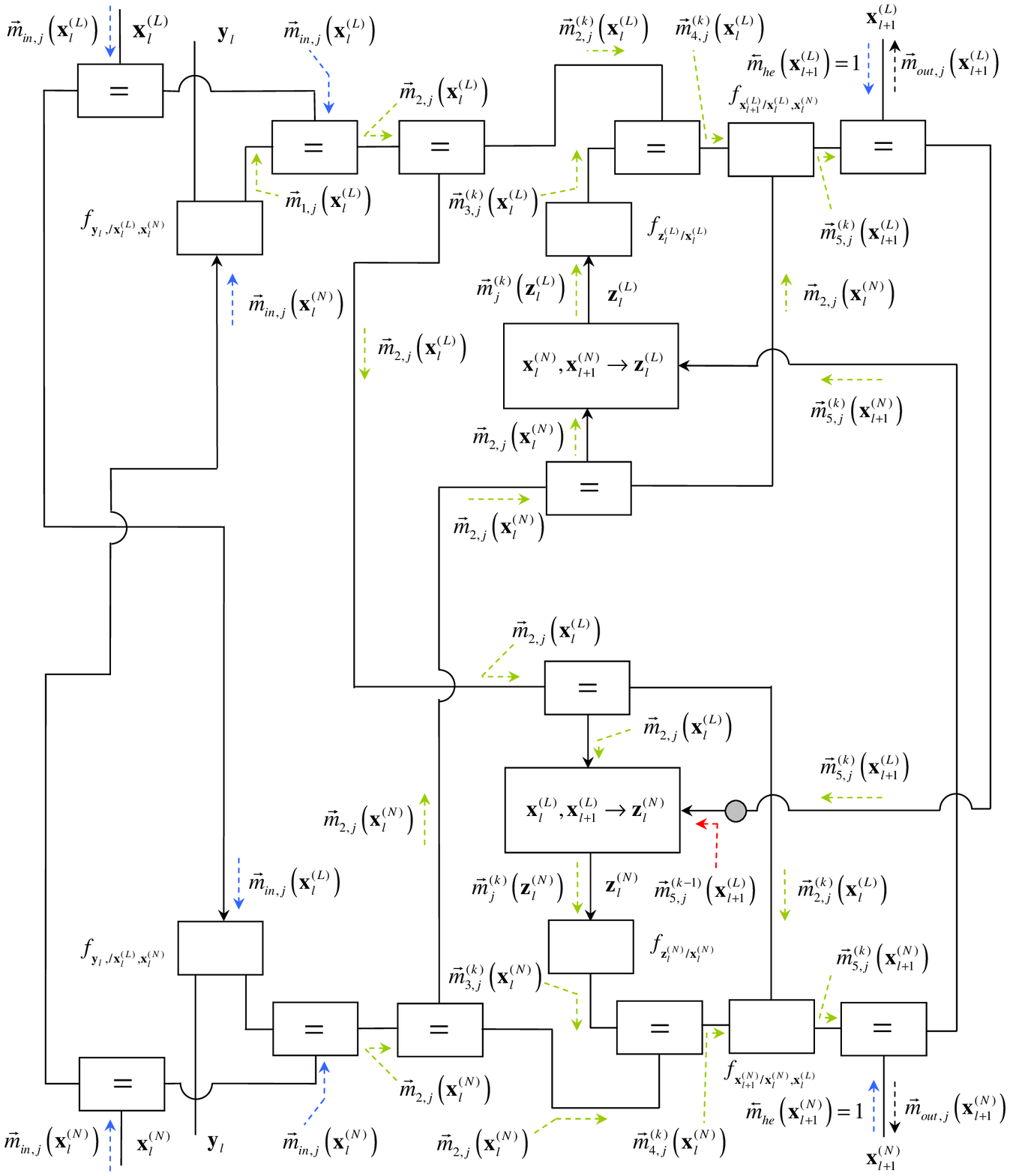}
	\caption{Message passing over the FG of Fig. 1 for the proposed TF technique. All the quantities appearing in this figure refer to the $k$-th iteration of
the $l$-th recursion; the messages available at the beginning of the
considered iteration are indicated by red arrows, those entering the graph by
blue arrows, those leaving it at the end of the last iteration by black arrows
and those computed within the considered iteration by green arrows.}
	\label{Fig_4}
\end{figure}

As already state above, our previous derivations rely on the assumption that
the message $\vec{m}_{5,j}(\mathbf{x}_{l+1}^{(L)})$ (\ref{eq:message_L_out})
is available when the time update for $\mathbf{x}_{l}^{(L)}$ is carried out;
unluckily, this message becomes available only in the last step of MPF.
However, if MPF is generalised in a way that $N_{it}>1$ iterations (i.e.,
message passes) are carried out within the same recursion, in the $k$-th
iteration (with $k=2,3,...,N_{it}$) the message $m_{j}(\mathbf{z}_{l}^{(N)})$
(\ref{eq:message_z_N1}) can be really evaluated exploiting the message
$\vec{m}_{5,j}(\mathbf{x}_{l+1}^{(L)})$ computed in the previous iteration.
These considerations lead, in a natural fashion, to the development of the
message passing illustrated in Fig. \ref{Fig_4}, which describes the message
flow occurring in the $k$-th iteration of a new filtering technique,
generalising MPF and called \emph{turbo filtering} (TF) in the following; note
that the superscripts $(k)$ or $(k-1)$ have been added to all the messages
flowing in the considered graph to identify the iteration in which they are
generated and that the grey circle appearing in the figure represents a unit
delay cell. The processing tasks accomplished by TF can be summarized as
follows. The first part of \ this technique can be considered as a form of
\emph{initialization}, in which the messages $\{\vec{m}_{2,j}(\mathbf{x}%
_{l}^{(N)})\}$ and $\{\vec{m}_{2,j}(\mathbf{x}_{l}^{(L)})\}$ are computed;
however, unlike MPF, resampling is not accomplished in the last part of the
time update for $\mathbf{x}_{l}^{(N)}$ (so that the messages $\{\vec{m}%
_{2,j}(\mathbf{x}_{l}^{(N)})\}$ are expressed by (\ref{eq:message_2_N_l-1})
instead of (\ref{eq:message_2_N_k-2}) and refer to the particle set
$S_{l/(l-1)}^{(N)}$). Moreover, as shown in Fig. \ref{Fig_4}, the messages
$\{\vec{m}_{2,j}(\mathbf{x}_{l}^{(N)})\}$ and $\{\vec{m}_{2,j}(\mathbf{x}%
_{l}^{(L)})\}$ emerging from the first part of TF remain unchanged within the
$l$-th recursion, since they represent the \emph{a priori} information
available about $\mathbf{x}_{l}^{(N)}$ and $\mathbf{x}_{l}^{(L)}$,
respectively; consequently, like in any turbo processing method, these
information are made available to all the iterations carried out within each
recursion. In the second part of TF\ $N_{it}$-iterations are accomplished with
the aim of progressively refining the set $S_{(l+1)/l}^{(N)}$ (the version of
this set generated in the $k$-th iteration is denoted $S_{(l+1)/l}^{(N)}[k]$)
and the associated Gaussian messages $\{\vec{m}_{5,j}^{(k)}(\mathbf{x}%
_{l+1}^{(L)})\}$. \ To achieve these result, in the $k $-th iteration (with
$k=1,2,...,N_{it}$) the ordered computation of the following messages is
accomplished: $\vec{m}_{j}^{(k)}(\mathbf{z}_{l}^{(N)})$ (\ref{eq:message_z_N1}%
), $\vec{m}_{3,j}^{(k)}(\mathbf{x}_{l}^{(N)})$ (\ref{eq:P_1_l_j-1}) (conveying
the weight $p_{l,j}[k]$), $\vec{m}_{4,j}^{(k)}(\mathbf{x}_{l}^{(N)})$
(\ref{eq:message_4_N_l-1-1}) (conveying the weight $W_{l,j}[k]$), $\vec
{m}_{5,j}^{(k)}(\mathbf{x}_{l+1}^{(N)})$ (\ref{eq:message_5_N}), $\vec{m}%
_{j}(\mathbf{z}_{l}^{(L)})$\ (\ref{eq:message_z_L}), $\vec{m}_{3,j}%
^{(k)}(\mathbf{x}_{l}^{(L)})$ (\ref{eq:message_3_L}), $\vec{m}_{4,j}%
^{(k)}(\mathbf{x}_{l}^{(L)})$ (\ref{eq:message_4-3}) and $\vec{m}_{5,j}%
^{(k)}(\mathbf{x}_{l+1}^{(L)})$ (\ref{eq:message_L_out}). Moreover, in the
$k$-th iteration resampling\footnote{Note that, after carrying out resampling,
the set of messages $\{\vec{m}_{2,j}(\mathbf{x}_{l}^{(L)})\}$ needs to be
properly reordered, since the messages associated with the discarded particles
are not preserved. This modifies the set of particles available in the next
iteration and, consequently, the set of weights $\{w_{l,j}\}$ associated with
them (these weights need to be renormalized after any change). In the
following the notation $\{w_{l,j}[k]\}$ is adopted to denote the set of
weights employed in the $k$-th iteration for the evaluation of the overall
weights $W_{l,j}[k] $ according to (\ref{eq:weight_before_resampling-1}).} is
accomplished on the basis of the particle weights $\{W_{l,j}[k]\}$; generally
speaking, this results in a new particle set denoted $S_{l/l}^{(N)}%
[k]=\{\mathbf{x}_{l/(l-1),j}^{(N)}[k],\,j=0,1,...,N_{p}-1\}$ (which is always
a subset of $S_{l/(l-1)}^{(N)}$). It is also important to point out that:

\begin{itemize}
\item In the first iteration (i.e, for $k=1$) the messages $\{\vec{m}%
_{5,j}^{(0)}(\mathbf{x}_{l+1}^{(N)})\}$ are undefined, so that $\{\vec{m}%
_{j}^{(1)}(\mathbf{z}_{l}^{(N)})=1\}$ (and, consequently, $\{\vec{m}%
_{3,j}^{(1)}(\mathbf{x}_{l}^{(N)})=1\}$) must be assumed; moreover, resampling
is not accomplished (i.e., $S_{l/l}^{(N)}[1]=S_{l/l}^{(N)}=S_{l/(l-1)}^{(N)}%
$), since the weights $p_{l,j}$ appearing in the overall weight $W_{l,j}$
(\ref{eq:weight_before_resampling-1}) become available in the following iterations.

\item In each iteration the equality $\overleftarrow{m}_{he}(\mathbf{x}%
_{l+1}^{(L)})=$ $\overleftarrow{m}_{he}(\mathbf{x}_{l+1}^{(N)})=1$ is assumed
for the two messages entering the FG along the half edges associated with
$\mathbf{x}_{l+1}^{(L)}$ and $\mathbf{x}_{l+1}^{(N)}$ (see Fig. \ref{Fig_4}),
since no information comes from the next recursion. For this reason, at the
end of the last iteration (i.e., for $k=N_{it}$), the output messages (i.e.,
the input messages feeding the $(l+1)$-th recursion) are evaluated as (see
(\ref{eq:message_5_N_l+1-1})-(\ref{eq:message_5_N_l+1}) and
(\ref{eq:message_L_out}))
\begin{equation}
\vec{m}_{out,j}\left(  \mathbf{x}_{l+1}^{(N)}\right)  =\vec{m}_{5,j}%
^{(N_{it})}\left(  \mathbf{x}_{l+1}^{(N)}\right)  \,\overleftarrow{m}%
_{he}(\mathbf{x}_{l+1}^{(N)})=\vec{m}_{5,j}^{(N_{it})}\left(  \mathbf{x}%
_{l+1}^{(N)}\right) \label{eq:message_5_N_l+1-1-1}%
\end{equation}
and
\begin{equation}
\vec{m}_{out,j}\left(  \mathbf{x}_{l+1}^{(L)}\right)  =\vec{m}_{5,j}%
^{(N_{it})}\left(  \mathbf{x}_{l+1}^{(L)}\right)  \,\overleftarrow{m}%
_{he}(\mathbf{x}_{l+1}^{(L)})=\vec{m}_{5,j}^{(N_{it})}\left(  \mathbf{x}%
_{l+1}^{(L)}\right)  ,\label{eq:message_L_out-2}%
\end{equation}
with $j=0,1,...,N_{p}-1$.
\end{itemize}%


Another relevant issue concerns the interpretation of the processing tasks
accomplished in the TF technique. In fact, our derivations show that, at the
end of the $k$-th iteration, the \emph{a posteriori} statistical information
about the $j$-th particle $\mathbf{x}_{l/l,j}^{(N)}[k]$ of $S_{l/l}^{(N)}[k]$
is provided by the message $\vec{m}_{4,j}^{(k)}(\mathbf{x}_{l}^{(N)})$
(\ref{eq:message_4_N_l-1-1}), which conveys the weight (see
(\ref{eq:weight_before_resampling-1}))%
\begin{equation}
W_{l,j}[k]=w_{l,j}[k]\cdot p_{l,j}[k]\cdot w_{l,j}^{(a)},\label{m_4}%
\end{equation}
where $w_{l,j}^{(a)}$ denotes the \emph{a priori} information available at the
beginning of the $l$-th recursion for the $j$-th particle ( in our derivation
$w_{l,j}^{(a)}$ $=1$ has been assumed, in place of $w_{l,j}^{(a)}$ $=1/N_{p}$,
to simplify the notation; see (\ref{eq:message_N_pred-1})), $p_{l,j}[k]$ is
the weight originating from $\vec{m}_{j}(\mathbf{z}_{l}^{(N)}) $
(\ref{eq:message_z_N1}) and conveyed by $\vec{m}_{3,j}^{(k)}(\mathbf{x}%
_{l}^{(N)})$ (\ref{eq:P_1_l_j-1}), and $w_{l,j}[k]$ is the weight computed on
the basis of the available measurement $\mathbf{y}_{l}$. Taking the natural
logarithm of both sides of (\ref{m_4}) produces
\begin{equation}
L_{l,j}[k]=L_{l,j}^{(a)}+L_{l,j}^{(y)}[k]+L_{l,j}^{(z)}[k],\label{loglik}%
\end{equation}
where $L_{l,j}[k]\triangleq\ln(W_{l,j}[k])$, $L_{l,j}^{(y)}[k]\triangleq
\ln(w_{l,j}[k])$, $L_{l,j}^{(z)}\triangleq\ln(p_{l,j}[k])$ and $L_{l,j}%
^{(a)}\triangleq\ln(w_{l,j}^{(a)})$. The last equation has exactly the same
structure as the well known formula (see \cite[Sec. 10.5, p. 450, eq.
(19.15)]{Vitetta} or \cite[Par. II.C, p. 432, eq. (20)]{Hagenauer_1996})
\begin{equation}
L\left(  u_{j}|\mathbf{y}\right)  =L\left(  u_{j}\right)  +L_{c}(y_{j}%
)+L_{e}\left(  u_{j}\right) \label{log_turbo}%
\end{equation}
expressing of the \emph{log-likelihood ratio} (LLR) available for the $j$-th
information bit \ $u_{j}$ at the output of a \emph{soft-input soft-output}
channel decoder operating\ over an \emph{additive white Gaussian noise} (AWGN)
channel and fed by: a) the channel output vector $\mathbf{y}$ (whose $j$-th
element $y_{j}$ is generated by the communications channel in response to a
channel symbol conveying $u_{j}$ and is processed to produce the so-called
\emph{channel} LLR $L_{c}(y_{j})$); b) the a priori LLR $L\left(
u_{j}\right)  $ about $u_{j}$; c) the \emph{extrinsic} LLR $L_{e}\left(
u_{j}\right)  $, i.e. a form of soft information available about $u_{j}$, but
intrinsically not influenced by such a bit (in turbo decoding of concatenated
channel codes extrinsic infomation is generated by another channel decoder
with which soft information is exchanged with the aim of progressively
refining data estimates). This correspondence is not only formal, since in
eqs. (\ref{loglik}) and (\ref{log_turbo}) terms playing similar roles can be
easily identified. For instance, the term $L_{l,j}^{(y)}[k]$ ($L_{l,j}^{(a)}$)
in (\ref{loglik}) provides the same kind of information as $L_{c}(y_{j})$
($L\left(  u_{j}\right)  $), since these are both related to the noisy data (a
priori information) available about the quantities to be estimated (the system
state in one case, an information bit the in the other one). What about the
term $L_{l,j}^{(z)}[k]$ appearing in the RHS of (\ref{loglik})? The link we
have established between (\ref{loglik}) and (\ref{log_turbo}) unavoidably
leads to the conclusion that such a term should represent the counterpart of
the quantity $L_{e}\left(  u_{j}\right)  $ appearing in (\ref{log_turbo}),
i.e. the so called \emph{extrinsic information} (in other words, that part of
the information\ available about $\mathbf{x}_{l}^{(N)}$ and not
\emph{intrinsically} influenced by $\mathbf{x}_{l}^{(N)}$ itself). This
interpretation is confirmed by the fact that $L_{l,j}^{(z)}[k]$ is computed on
the basis of the statistical knowledge available about $\mathbf{x}_{l}^{(L)}$
and $\mathbf{x}_{l+1}^{(L)}$ (see (\ref{eq:av_z}) and (\ref{eq:C_z-1})),
which, thanks to (\ref{eq:XL_update}), does provide useful information about
$\mathbf{x}_{l}^{(N)}$. The theory of turbo decoding of channel codes shows
that, generally speaking, extrinsic information originates from code
constraints. In our scenario, a similar interpretation can be also provided
for $L_{l,j}^{(z)}[k]$, since $\mathbf{x}_{l+1}^{(L)}$ can be seen as the
noisy output of a communication channel, affected by the bias $\mathbf{f}%
_{l,j}^{(L)}$ and the additive noise $\mathbf{w}_{l}^{(L)}$, and over which
the codeword $\mathbf{A}_{l,j}^{(L)}\mathbf{x}_{l}^{(L)}$ of a rate-$1$ block
code is transmitted in response to the message $\mathbf{x}_{l}^{(L)}$. These
considerations show that, in evaluating $L_{l,j}^{(z)}[k]$, we are actually
exploiting a sort of `code' constraints, which are mathematically expressed by
(\ref{eq:XL_update}). The reader can easily verify that a similar
interpretation can be provided for $\vec{m}_{j}(\mathbf{z}_{l}^{(L)})$
(\ref{eq:message_z_L}), which represents the \emph{extrinsic information}
component\footnote{In practice, the mechanism employed to generate
$\mathbf{z}_{l,j}^{(L)}$ is based on the state update equation
(\ref{eq:XN_update}).} contained in $\vec{m}_{4,j}(\mathbf{x}_{l}^{(L)})$
(\ref{eq:message_4-3}) (conveying our \emph{a posteriori} information about
$\mathbf{x}_{l}^{(L)}$); the other two components are represented by the
message $\vec{m}_{2,j}(\mathbf{x}_{l}^{(L)})$ (\ref{eq:message_2-l_L})
(representing the \emph{measurement information} about $\mathbf{x}_{l}^{(L)}$)
and the message $\vec{m}_{in,j}(\mathbf{x}_{l}^{(L)})$
(\ref{eq:eq:message_L_pred-1}) (corresponding to our\emph{\ a priori}
information about $\mathbf{x}_{l}^{(L)}$). Consequently, TF can be seen, in
the domain of Bayesian filtering techniques, as the counterpart of turbo
decoding of concatenated codes; this parallelism can be exploited to provide
further insights into iterative filtering techniques. For instance, it is well
known that, in turbo decoding of concatenated channel codes, the extrinsic
information generated by soft decoders become more and more correlated as
iterations evolve;\ this entails that diminishing benefits are provided by
additional iterations. This phenomenon should be observed in TF too for
similar reasons and can be motivated by rewriting $\mathbf{\eta}%
_{\mathbf{z},l,j}^{(N)}$ (\ref{eq:av_z}) and $\mathbf{C}_{\mathbf{z}%
,l,j}^{(N)}$ (\ref{eq:C_z-1}) as%

\begin{equation}
\mathbf{\eta}_{\mathbf{z},l,j}^{(N)}=\mathbf{f}_{l,j}^{(L)}+\mathbf{A}%
_{l,j}^{(L)}\left[  \mathbf{\eta}_{4,l,j}^{(L)}-\mathbf{\eta}_{2,l,j}%
^{(L)}\right] \label{eq:av_z_final}%
\end{equation}
and%
\begin{equation}
\mathbf{C}_{\mathbf{z},l,j}^{(N)}=\mathbf{C}_{w}^{(L)}+\mathbf{A}_{l,j}%
^{(L)}\left[  \mathbf{C}_{4,l,j}^{(L)}-\mathbf{C}_{2,l,j}^{(L)}\right]
\left(  \mathbf{A}_{l,j}^{(L)}\right)  ^{T},\label{eq:C_Z_final}%
\end{equation}
respectively (thanks to (\ref{eq:eta_5_new}) and (\ref{eq:C_5_new}),
respectively). In fact, the last two equations show that the vector
$\mathbf{\eta}_{\mathbf{z},l,j}^{(N)}$ and the matrix $\mathbf{C}%
_{\mathbf{z},l,j}^{(N)}$ are influenced by the \emph{difference} between the
statistical information (expressed by a mean vector and a covariance matrix)
available about $\mathbf{x}_{l}^{(L)}$ before processing the pseudomeasurement
$\mathbf{z}_{l}^{(L)}$ and those available after this task has been carried
out. In other words, the extrinsic information provided by $\mathbf{z}%
_{l}^{(L)}$ influences $\mathbf{z}_{l}^{(N)}$ and viceversa.

Finally, it is important to point out that the proposed analogy between turbo
filtering and turbo decoding suggests the potential limits of the TF technique
(and of any other iterative filtering method relying on the developed FG). In
fact, it is well known that turbo decoding methods do not provide real
benefits below a certain signal-to-noise ratio, i.e. when the quality of the
received signal is so poor that the transmitted coded sequence cannot be
recovered. A similar phenomenon is expected occur with TF too; consequently,
this filtering method could not outperform MPF in the presence of strong
measurement noise and/or fast dynamics affecting the considered SSM.

\section{Numerical Results\label{num_results}}

In this Section MPF and the related filtering methods developed in this
manuscript are compared in terms of accuracy and computational load for a
specific CLG system, characterized by $D_{L}=3$, $D_{N}=1$ (so that $D=4$) and
$P=2$. The structure of the considered system has been partly inspired by the
example proposed by Sch\"{o}n in \cite{Schon_2010} and is characterized by: a)
the state models
\begin{equation}
\mathbf{x}_{l+1}^{(L)}=\left(
\begin{array}
[c]{ccc}%
0.8 & 0.2 & 0\\
0 & 0.7 & -0.2\\
0 & 0.2 & 0.7
\end{array}
\right)  \mathbf{x}_{l}^{(L)}+\left(
\begin{array}
[c]{c}%
\cos(x_{l}^{(N)})\\
-\sin(x_{l}^{(N)})\\
0.5\sin(2x_{l}^{(N)})
\end{array}
\right)  +\mathbf{w}_{l}^{(L)}\label{state_mod_1}%
\end{equation}
and%
\begin{equation}
x_{l+1}^{(N)}=\arctan\left(  x_{l}^{(N)}\right)  +\left(  0.9\quad
0\quad0\right)  \mathbf{x}_{l}^{(L)}+w_{l}^{(N)}\label{state_mod_2}%
\end{equation}
with $\mathbf{w}_{l}^{(L)}\sim\mathcal{N}(0,(\sigma_{w}^{(L)})^{2}%
\mathbf{I}_{3})$, $w_{l}^{(N)}\sim\mathcal{N}(0,(\sigma_{w}^{(N)})^{2}$; b)
the measurement model%
\begin{equation}
\mathbf{y}_{l}=\left(
\begin{array}
[c]{c}%
0.1\left(  x_{l}^{(N)}\right)  {}^{2}\cdot\text{sgn}\left(  x_{l}^{(N)}\right)
\\
0
\end{array}
\right)  +\left(
\begin{array}
[c]{ccc}%
0 & 0 & 0\\
1 & -1 & 1
\end{array}
\right)  \mathbf{x}_{l}^{(L)}+\mathbf{e}_{l}\label{meas_syst}%
\end{equation}
with $\mathbf{e}_{l}\sim\mathcal{N}(0,(\sigma_{e})^{2}\mathbf{I}_{2})$. Note
that the state equation (\ref{state_mod_1}), unlike its counterpart proposed
in \cite{Schon_2010}, depends on $x_{l}^{(N)}$ (i.e., it contains a function
$\mathbf{f}_{l}^{(L)}\left(  x_{l+1}^{(N)}\right)  \neq\mathbf{0}_{3} $ in its
RHS), so that TF, which relies on the availability of the vector
$\mathbf{z}_{l}^{(N)}$ (\ref{eq:z_N_l}), can be employed for this system.

In our computer simulations the \emph{root mean square error} (RMSE) has been
evaluated to compare the accuracy of the state estimates generated by
different filtering techniques. More specifically, for each technique two
RMSEs have been computed, one (denoted $RMSE_{L}($alg$)$, where `alg' denotes
the algorithm this parameter refers to) representing the square root of the
average \emph{mean square error} (MSE) evaluated for the three elements of
$\mathbf{x}_{l}^{(L)}$, the other one (denoted $RMSE_{N}($alg$)$) referring to
the (monodimensional) nonlinear component $x_{l}^{(N)}$; this distinction is
important since, as shown by our simulation results, the estimation accuracy
for $\mathbf{x}_{l}^{(L)}$ can be quite different from (and is usually smaller
than) that referring to $x_{l}^{(N)}$.

As far the assessment of the computational requirements of the investigated
filtering techniques is concerned, MPF (for which an accurate analysis of its
computational complexity is available in \cite{Schon_2005_complexity}) has
been taken as a \emph{baseline}. For this reason, our comparisons between the
considered filtering techniques are based on the evaluation of a single
parameter, denoted $\Delta_{c}($alg$)$ and representing the \emph{percentage
variation in computation time of the considered algorithm }(denoted
`alg')\emph{\ with respect to} MPF (operating with the same parameters and, in
particular, with the same $N_{p}$ as alg).

Moreover, in our computer simulations, the following choices have been made:
a) $\sigma_{w}^{(L)}=\sigma_{w}^{(N)}=5\cdot10^{-3}$ has been selected for the
standard deviations of the process noises $\{\mathbf{w}_{k}^{(L)}\}$ and
$\{w_{k}^{(N)}\}$, unless differently stated; b) the so called \emph{jittering
technique} \cite{Li_2015} has been employed to mitigate the so called
\emph{depletion problem} in the generation of new particles.

Some results illustrating the dependence of $RMSE_{L}$ and $RMSE_{N}$ on the
number of particles ($N_{p}$) for MPF and SMPF\#1 are illustrated in Fig.
\ref{Fig_1_sim} (in this and in the following figures simulation results are
identified by markers, whereas continuous lines are drawn to ease reading);
$\sigma_{e}=10^{-2}$ and $N_{p}\in\lbrack50,300]$ have been selected in this
case. From these results the following conclusions can be easily inferred for
the considered system:

\begin{enumerate}
\item A negligible improvement in the estimation accuracy of both MPF and
SMPF\#1 is achieved if the value of $N_{p}$ exceeds $200$.

\item A significant gap between $RMSE_{L}($MPF$)$ and $RMSE_{N}($MPF$)$
exists; this is motivated by the fact that, in MPF, the estimation of
$\mathbf{x}_{l}^{(L)}$ relies on both the real measurement $\mathbf{y}_{l}$
(\ref{meas_syst}) and the pseudo-measurements $\{\mathbf{z}_{l,j}^{(L)}\}$
(\ref{eq:message_Z_L}), whereas the estimation of $x_{l}^{(N)}$ benefits from
$\mathbf{y}_{l}$ only.

\item The gap between $RMSE_{L}($SMPF\#1$)$ and $RMSE_{N}($SMPF\#1$)$ is much
smaller than that observed for MPF, even if the pseudo-measurements
$\{\mathbf{z}_{l,j}^{(L)}\}$ are also exploited by SMPF\#1. \ This reduction
in the RMSE gap can be related to the degradation in the estimation accuracy
of $\mathbf{x}_{l}^{(L)}$; for instance, for $N_{p}=200$, $RMSE_{L}%
($SMPF\#1$)$ is about twice $RMSE_{L}($MPF$)$ (on the contrary, $RMSE_{N}%
($SMPF\#1$)\cong1.09\cdot$ $RMSE_{N}($MPF$)$).
\end{enumerate}

Our numerical results have also evidenced that: a) SMPF\#1 requires a
substantially smaller computational effort than MPF, since $\Delta_{c}%
($SMPF\#1$)$ approximately ranges in the interval $[-62\%,-58\%]\ $ for the
considered values of $N_{p}$; b) despite the above mentioned degradation in
estimation accuracy (see point 3.), SMPF\#1 does not suffer from
\emph{tracking losses} in the considered scenario; c) SMPF\#2 accuracy is
almost indentical that of SMPF\#1; d) $\Delta_{c}($SMPF\#2$)$ is approximately
lower than $\Delta_{c}($SMPF\#1$)$ by $6\%$ in the considered range for
$N_{p}$ and, consequently, achieves a better complexity-performance tradeoff
than SMPF\#1. All this suggests that simplified MPF techniques can be really
developed without incurring the serious technical problems that affect the
techniques proposed in \cite{Mustiere_2006} (tracking losses and poor RMSE
performance) and in \cite{Lu_2007} (partitioning and update of the particle
set into groups within each recursion; see Section \ref{sec:simplifying}).%
\begin{figure}
	\centering
		\includegraphics{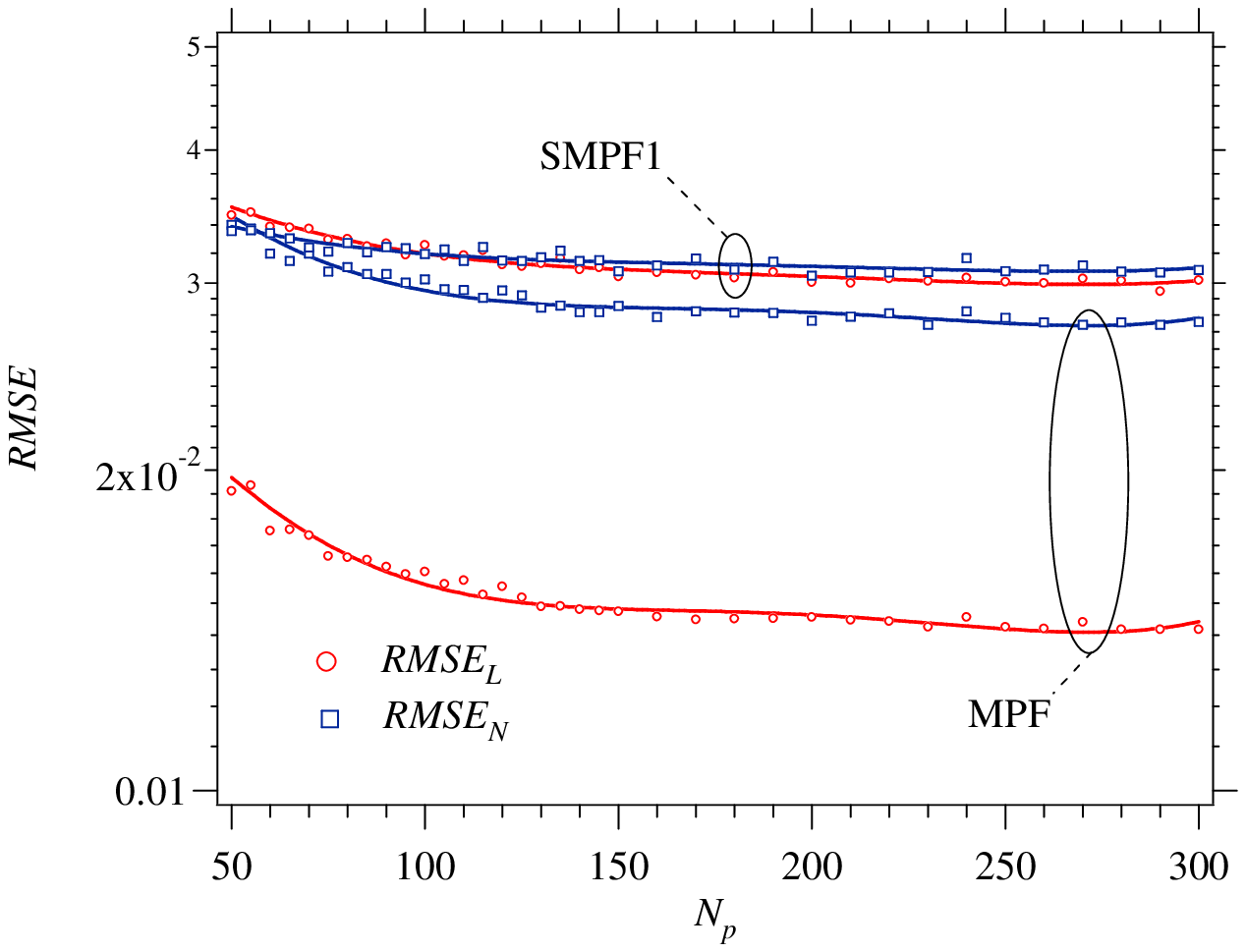}
	\caption{RMSE performance versus $N_{p}$ for the linear component ($RMSE_{L}$)
and the nonlinear component ($RMSE_{N}$) of the state $\mathbf{x}_{l}$ for the
system described by eqs. (\ref{state_mod_1})-(\ref{meas_syst}). MPF and
SMPF\#1 are considered; in both cases $\sigma_{w}^{(L)}=\sigma_{w}%
^{(N)}=5\cdot10^{-3}$ and $\sigma_{e}=10^{-2}$ have been selected.}
	\label{Fig_1_sim}
\end{figure}

\begin{figure}
	\centering
		\includegraphics{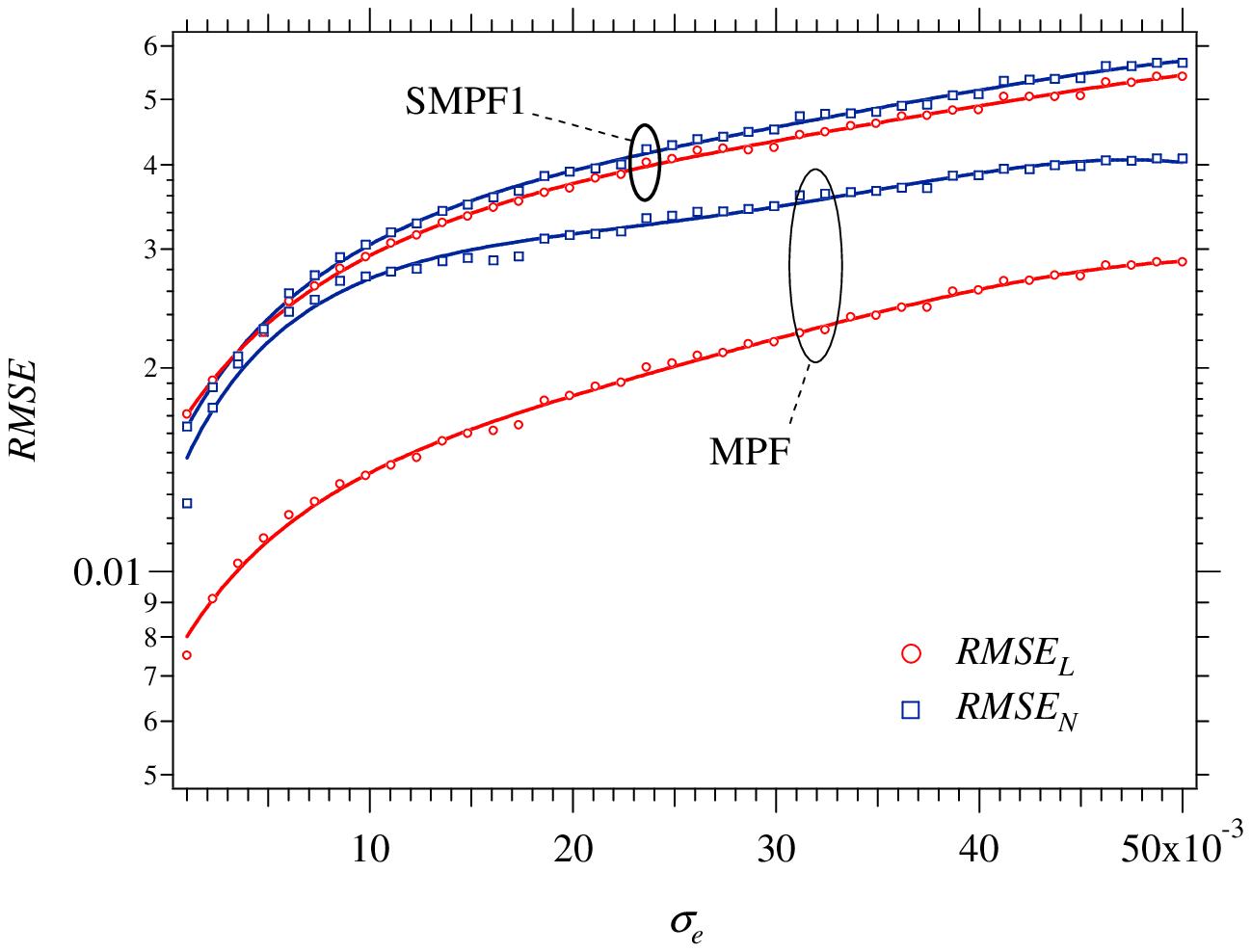}
	\caption{RMSE performance versus $\sigma_{e}$ for the linear component
($RMSE_{L}$) and the nonlinear component ($RMSE_{N}$) of the state
$\mathbf{x}_{l}$ for the system described by eqs. (\ref{state_mod_1}%
)-(\ref{meas_syst}). MPF and SMPF\#1 are considered; in both cases $\sigma
_{w}^{(L)}=\sigma_{w}^{(N)}=5\cdot10^{-3}$ and $N_{p}=200$ have been
selected.}
	\label{Fig_2_sim}
\end{figure}

The dependence of $RMSE_{L}$ and $RMSE_{N}$ on $\sigma_{e}$ (i.e., on the
intensity of the noise affecting the available measurements) has been also
analysed for MPF and SMPF\#1. Some results are shown in Fig. \ref{Fig_2_sim};
$N_{p}$ $=500$ and $\sigma_{e}\in\lbrack10^{-3},5\cdot10^{-2}]$ have been
selected in this case. These results show that the gap between $RMSE_{L}%
($MPF$)$ and $RMSE_{L}($SMPF\#1$)$ slightly increases as $\sigma_{e} $ becomes
smaller; the opposite occurs for $RMSE_{N}($MPF$)$ and $RMSE_{N}($SMPF\#1$)$.
These results can be related again to the fact that the pseudo-measurements
$\{\mathbf{z}_{l,j}^{(L)}\}$ (\ref{eq:message_Z_L}) really play a role in the
MPF estimation of $\mathbf{x}_{l}^{(L)}$ and that the quality of the
information conveyed by these pseudo-measurements indirectly improves as the
real measurements become less noisy.%

\begin{figure}
	\centering
		\includegraphics{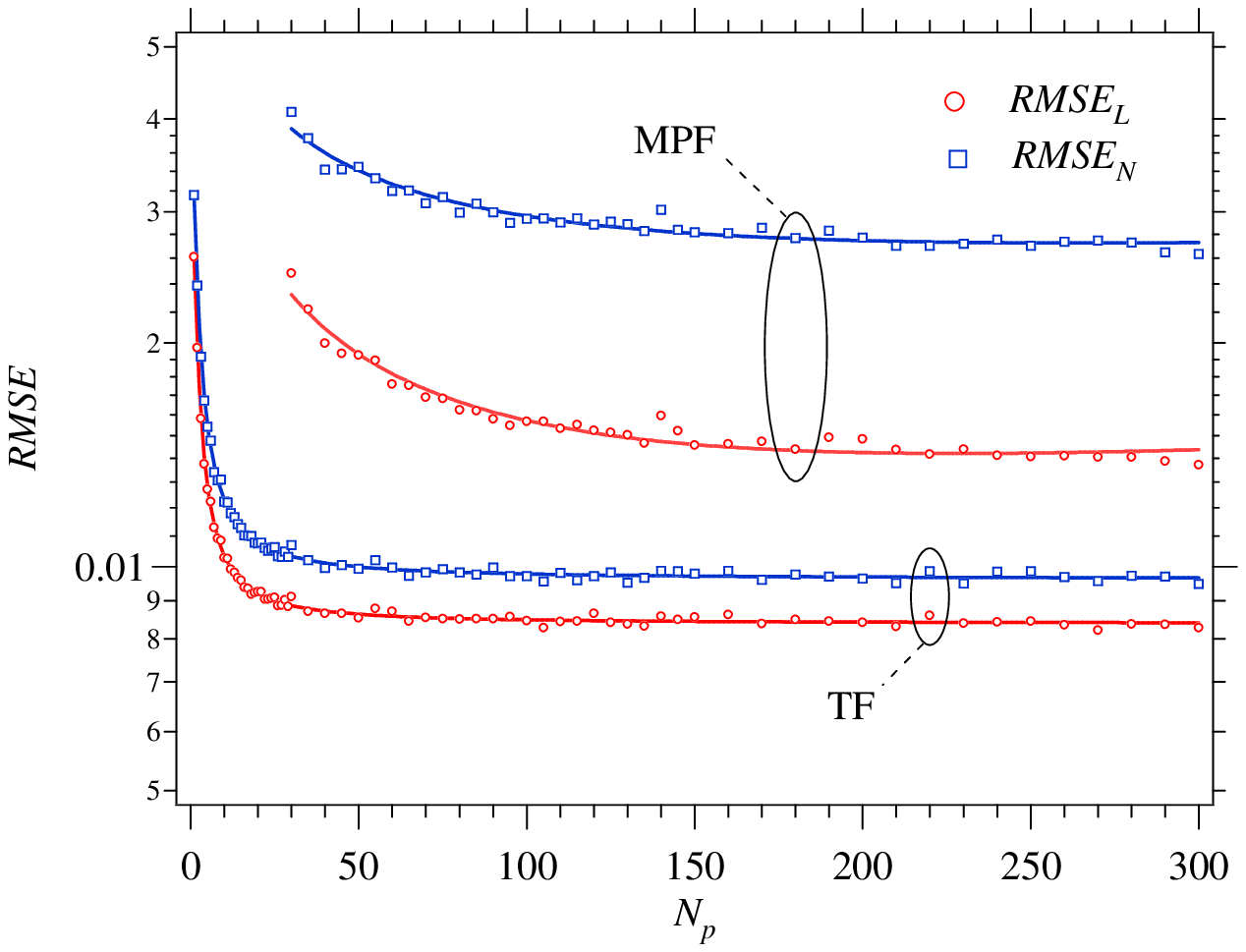}
	\caption{RMSE performance versus $N_{p}$ for the linear component ($RMSE_{L}$)
and the nonlinear component ($RMSE_{N}$) of the state $\mathbf{x}_{l}$ for the
system described by eqs. (\ref{state_mod_1})-(\ref{meas_syst}). MPF and TF are
considered; in both cases $\sigma_{w}^{(L)}=\sigma_{w}^{(N)}=5\cdot10^{-3}$
and $\sigma_{e}=10^{-2}$ have been selected.}
	\label{Fig_3sim}
\end{figure}

Some results illustrating the dependence of $RMSE_{L}$ and $RMSE_{N}$ on the
number of particles ($N_{p}$) for MPF and TF are illustrated in Fig.
\ref{Fig_3sim}; $\sigma_{e}=10^{-2}$ , $N_{p}\in\lbrack1,300]$ and $N_{it}=2$
for TF, and $N_{p}\in\lbrack30,300]$ for MPF have been chosen in this case.
From these results it is easily inferred that:

\begin{enumerate}
\item TF outperforms MPF in tems of both $RMSE_{L}$ and $RMSE_{N}$; for
instance, $RMSE_{L}($MPF$)\cong1.71\cdot$ $RMSE_{L}($TF$)$ and $RMSE_{N}%
($MPF$)\cong2.86\cdot$ $RMSE_{N}($TF$)$ if $N_{p}=200$ is selected.

\item The gap between $RMSE_{L}($MPF$)$ and $RMSE_{N}($MPF$)$ is substantially
larger than the corresponding gap for TF (in particular, $RMSE_{N}%
($MPF$)-RMSE_{L}($MPF$)$ $\cong13.4\cdot\lbrack RMSE_{N}($TF$)-RMSE_{L}%
($TF$)]$; this can be easily related to the fact that, unlike MPF, in TF the
estimation of both $\mathbf{x}_{l}^{(L)}$ and $x_{l}^{(N)}$ benefits from the
availability of pseudo-measurements.

\item The performance gap between TF and MPF increases as $N_{p}$ gets
smaller. Moreover, TF accuracy starts quickly degrading when $N_{p}$ drops
below $11$, whereas the same phenomenon starts for $N_{p}\cong30$ with MPF.
Note, however, that different phenomena occur with MPF and TF when the
particle set is small. In fact, our computer simulations have shown that MPF
suffers from frequent tracking losses for $N_{p}<30$ (since it is unable to
generate a reliable representation of $x_{l}^{(N)}$ when a limited number of
particles is available); on the contrary, TF does not suffer from the same
problem even if a very small particle set is used.
\end{enumerate}

We believe that last result is really important and can be motivated as
follows. The TF technique, through its feedback mechanism, makes a
substantially more efficient use of the available particles than MPF; this
results in an appreciable improvement of both stability and accuracy of state estimation.%

Our simulation results have also shown that: a) in the considered scenario a
negligible improvement in estimation accuracy is obtained if $N_{it}>2$ is
selected; b) $\Delta_{c}($TF$)$ approximately ranges in the interval
$[70\%,80\%]\ $ for the considered values of $N_{p}$ and, consequently,
requires a substantially larger computational effort than MPF if these
algorithms operate with the same number of particles. In practice, however, as
evidenced by the results shown in Fig. \ref{Fig_3sim}, TF can reliably operate
with a very small particle set and, consequently, \emph{it outperforms }MPF
\emph{in terms of both performance and complexity} if $N_{p} $ is properly
selected. For instance, in the considered scenario TF with $N_{p}=20$
particles achieves a substally better accuracy than MPF with $N_{p}=40$
particles, even if, as evidenced by our computer simulations, they
approximately require the same computation time. It is not difficult to show
that similar considerations hold if SMPF\#1 and SMPF\#2 are considered in
place of MPF. Therefore, our results suggest that the real key to the
complexity reduction of the filtering techniques relying on the FG\ shown in
Fig. \ref{Fig_2} is not provided by the approximations adopted for the MPF
processing tasks in Section \ref{sec:simplifying}, but by the exploitation of
the new pseudo-measurement $\mathbf{z}_{l}^{(N)}$ (\ref{eq:z_N_l}). We should
never forget, however, that TF cannot be adopted for all CLG systems; in fact,
it cannot be employed if $\mathbf{f}_{l}^{(L)}\left(  x_{l+1}^{(N)}\right)
=\mathbf{0}_{D_{L}}$ in the RHS of (\ref{eq:XL_update}).

\begin{figure}
	\centering
		\includegraphics{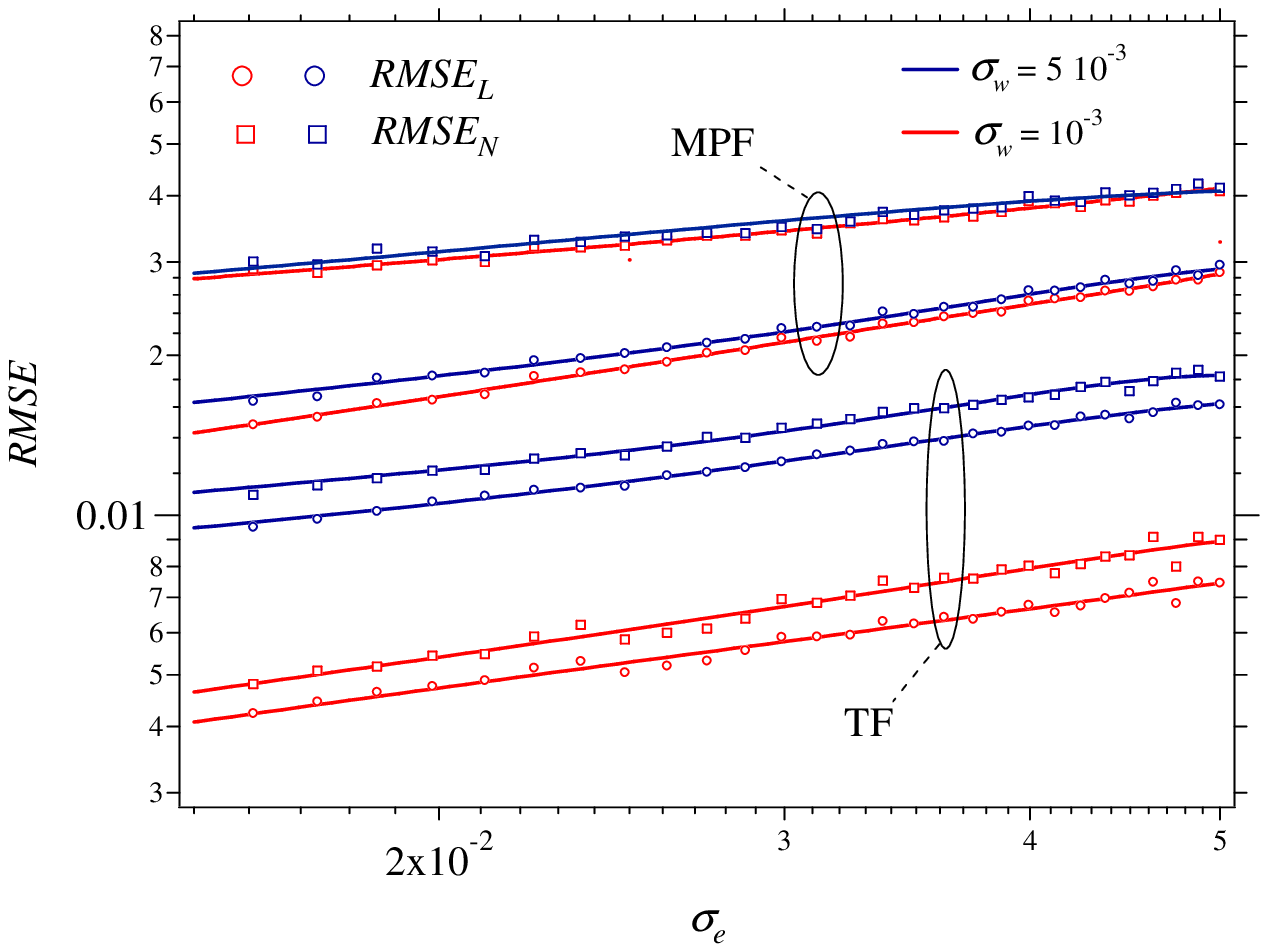}
	\caption{RMSE performance versus $\sigma_{e}$ for the linear component
($RMSE_{L}$) and the nonlinear component ($RMSE_{N}$) of the state
$\mathbf{x}_{l}$ for the system described by eqs. (\ref{state_mod_1}%
)-(\ref{meas_syst}). MPF and TF are considered; $\sigma_{e}=10^{-2}$ and
$N_{p}=200$ are assumed in both cases. As far system noise is concerned, the
cases $\sigma_{w}^{(L)}=\sigma_{w}^{(N)}=5\cdot10^{-3}$ and $\sigma_{w}%
^{(L)}=\sigma_{w}^{(N)}=10^{-1}$ are taken into consideration.}
	\label{Fig_4_sim}
\end{figure}

Finally, the dependence of $RMSE_{L}$ and $RMSE_{N}$ on $\sigma_{e}$ has been
assessed for TF and compared with that characterizing MPF. Some numerical are
illustrated in Fig. \ref{Fig_4_sim}. In this case, we have selected
$\sigma_{e}\in\lbrack1.5\cdot10^{-2},5\cdot10^{-2}]$ , $N_{p}=200$ and
$N_{it}=2$ for TF, and $N_{p}=200$ for MPF; moreover, $\sigma_{w}^{(L)}%
=\sigma_{w}^{(N)}=5\cdot10^{-3}$ and $\sigma_{w}^{(L)}=\sigma_{w}%
^{(N)}=10^{-3}$ have been considered to analyse the dependence of the
performance gap between TF and MPF on the intensity of process noise (and,
consequently, on system dynamics). These results show that: a) the performance
gap between TF and MPF undergoes small changes if $\sigma_{e}$ varies in the
considered interval; b) on the contrary, a substantial change in this gap is
obtained if $\sigma_{w}^{(L)}$ and $\sigma_{w}^{(N)}$ are reduced from
$5\cdot10^{-3}$ to $10^{-3}$. This suggests that measurement noise and process
noise can have different impacts on the performance gain provided by TF over
MPF.

\section{Conclusions\label{sec:conc}}

In this manuscript a FG approach has been employed to analyse the filtering
problem for mixed linear/nolinear models. This has allowed us to: a) prove
that this problem involves a FG which is not cycle free; b) provide a new
interpretation of MPF as a \emph{forward only} message passing algorithm over
a specific FG; c) develop novel filtering algorithms for simplifying or
generalising it. In particular, an important iterative filtering technique,
dubbed \emph{turbo filtering}, has been devised and its relation with the
turbo decoding techniques for concatenated channel codes has been analysed in
detail. All the considered filtering techniques have been compared in terms of
both accuracy and computational requirements for a specific CLG system. The
most interesting result emerging from our computer simulations is represented
by the clear superiority of turbo filtering over marginalized particle
filtering. In fact, the former technique, through the exploitation of new
pseudo-measurements, can achieve a better accuracy than the latter one at an
appreciably smaller computational load. Our ongoing research activities in
this area include the development of other related filtering techniques and
the application of turbo filtering to specific state estimation problems.

\appendix

\section{Appendix}

Given the pdfs $f_{1}(\mathbf{y})\triangleq\mathcal{N}(\mathbf{y}%
;\mathbf{\eta}_{1},\mathbf{C}_{1})$ and $f_{2}(\mathbf{y})\triangleq
\mathcal{N}(\mathbf{y};\mathbf{\eta}_{1},\mathbf{C}_{2})$ for the
$N$-dimensional vector $\mathbf{y}$, we are interested in evaluating the
\emph{correlation} between these two functions, i.e. the quantity%
\begin{equation}
c_{1,2}\triangleq\int f_{1}(\mathbf{y})\cdot f_{2}(\mathbf{y})\,d\mathbf{y}%
\label{eq:app_eq_0}%
\end{equation}
Substituting the expressions of $f_{1}(\mathbf{y})$ and $f_{2}(\mathbf{y})$ in
the RHS of the last equation produces, after some manipulation,
\begin{equation}
c_{1,2}=D\,\exp\left[  \frac{1}{2}\left(  \mathbf{\eta}^{T}\mathbf{W}%
\mathbf{\eta}-\mathbf{\eta}_{1}^{T}\mathbf{W}_{1}\mathbf{\eta}_{1}%
-\mathbf{\eta}_{2}^{T}\mathbf{W}_{2}\mathbf{\eta}_{2}\right)  \right]
\label{eq:app_eq_1}%
\end{equation}
where $\mathbf{\mathbf{W}}_{1}\triangleq\mathbf{C}_{1}^{-1}$,
$\mathbf{\mathbf{W}}_{2}\triangleq\mathbf{C}_{2}^{-1}$,
\begin{equation}
\mathbf{W}=\mathbf{\mathbf{W}}_{1}+\mathbf{\mathbf{W}}_{2},\label{eq:app_eq_2}%
\end{equation}%
\begin{equation}
\mathbf{W}\mathbf{\eta}=\mathbf{W}_{1}\mathbf{\eta}_{1}+\mathbf{W}%
_{2}\mathbf{\eta}_{2}\label{eq:app_eq_3}%
\end{equation}
and
\begin{equation}
D=\left(  2\pi\det\left[  \mathbf{C}_{1}+\mathbf{C}_{2}\right]  \right)
^{-N/2}.\label{eq:app_eq_4}%
\end{equation}


\begin{thebibliography}{99}                                                                                               %
\bibitem {Arulampalam_2002}M. S. Arulampalam, S. Maskell, N. Gordon and T.
Clapp, ``A Tutorial on Particle Filters for Online Nonlinear/Non-Gaussian
Bayesian Tracking'', \emph{IEEE Trans. Sig. Proc.}, vol. 50, no. 2, pp.
174-188, Feb. 2002.

\bibitem {Daum_1986}F. E. Daum, ``Exact Finite-Dimensional Nonlinear
Filters'', \emph{IEEE Tran. Aut. Contr.}, vol. 31, no. 7, pp. 616-622, July 1986.

\bibitem {Mazuelas_2013}S. Mazuelas, Y. Shen and M. Z. Win, ``Belief
Condensation Filtering'', \emph{IEEE Trans. Sig. Proc.}, vol. 61, no. 18, pp.
4403-4415, Sept. 2013.

\bibitem {Smidl_2008}V. Smidl and A. Quinn, \textquotedblleft Variational
Bayesian Filtering\textquotedblright, \emph{IEEE Trans. Sig. Proc.}, vol. 56,
no. 10, pp. 5020-5030, Oct. 2008.

\bibitem {Simandla_2006}M. \u{S}imandl, J. Kr\'{a}loveca and T.
S\"{o}derstr\"{o}mc, \textquotedblleft Advanced Point-Mass Method for
Nonlinear State Estimation\textquotedblright, \emph{Automatica}, vol. 42, pp.
1133-1145, 2006.

\bibitem {Anderson_1979}B. Anderson and J. Moore, \textbf{Optimal Filtering},
Englewood Cliffs, NJ, Prentice-Hall, 1979.

\bibitem {Julier_2004}S. J. Julier and J. K. Uhlmann, ``Unscented Filtering
and Nonlinear Estimation'', \emph{IEEE Proc.}, vol. 92, no. 3, pp. 401-422,
Mar. 2004.

\bibitem {Doucet_2001}A. Doucet, J. F. G. de Freitas and N. J. Gordon, ``An
Introduction to Sequential Monte Carlo methods,'' in \textbf{Sequential Monte
Carlo Methods in Practice}, A. Doucet, J. F. G. de Freitas, and N. J. Gordon,
Eds. New York: Springer-Verlag, 2001.

\bibitem {Doucet_2000}A. Doucet, S. Godsill and C. Andrieu, \textquotedblleft
On Sequential Monte Carlo Sampling Methods for Bayesian Filtering'',
\emph{Statist. Comput.}, vol. 10, no. 3, pp. 197-208, 2000.

\bibitem {Gustafsson_2010}F. Gustafsson, \textquotedblleft Particle Filter
Theory and Practice with Positioning Applications\textquotedblright,
\emph{IEEE Aerosp. and Electr. Syst. Mag.}, vol. 25, no. 7, pp. 53-82, July 2010.

\bibitem {Gustafsson_et_al_2002}F. Gustafsson, F. Gunnarsson, N. Bergman, U.
Forssell, J. Jansson, R. Karlsson and P. Nordlund, ``Particle Filters for
Positioning, Navigation, and Tracking'', \emph{IEEE Trans. Sig. Proc.}, vol.
50, 425-435, 2002.

\bibitem {Nordlund_2009}P. J. Nordlund and F. Gustafsson, \textquotedblleft
Marginalized Particle Filter for Accurate and Reliable Terrain-Aided
Navigation\textquotedblright, \emph{IEEE Trans. on Aerosp. and Elec. Syst.},
vol. 45, no. 4, pp. 1385-1399, Oct. 2009.

\bibitem {Bucy_1971}R. Bucy and K. Senne, \textquotedblleft Digital Synthesis
of Non-Linear Filters\textquotedblright, \emph{Automatica}, vol. 7, no. 3, pp.
287-298, 1971.

\bibitem {Daum_2003}F. Daum and J. Huang, \textquotedblleft Curse of
Dimensionality and Particle Filters\textquotedblright, \emph{Proc. IEEE
Aerosp. Conf.}, vol. 4, pp. 1979-1993, March 2003.

\bibitem {Schon_2005}T. Sch\"{o}n, F. Gustafsson, P.-J. Nordlund,
\textquotedblleft Marginalized Particle Filters for Mixed Linear/Nonlinear
State-Space Models\textquotedblright, \emph{IEEE Trans. Sig. Proc.}, vol. 53,
no. 7, pp. 2279-2289, July 2005.

\bibitem {Schon_2005_complexity}R. Karlsson, T. Sch\"{o}n, F. Gustafsson,
\textquotedblleft Complexity Analysis of the Marginalized Particle
Filter\textquotedblright, \emph{IEEE Trans. Sig. Proc.}, vol. 53, no. 11, pp.
4408-4411, Nov. 2005.

\bibitem {Mustiere_2006}F. Mustiere, M. Bolic and M. Bouchard, ``A Modified
Rao-Blackwellised Particle Filter'', \emph{Proc. of the 2006 IEEE Int. Conf.
on Ac., Sp. and Sig. Proc.} (ICASSP 2006), vol. 3, 14-19 May 2006.

\bibitem {Lu_2007}T. Lu, M. F. Bugallo and P. M. Djuric, \textquotedblleft
Simplified Marginalized Particle Filtering for Tracking Multimodal
Posteriors\textquotedblright, \emph{Proc. IEEE/SP 14th Workshop on Stat. Sig.
Proc. }(SSP '07), pp. 269-273, Madison, WI (USA), 2007.

\bibitem {Lindsten_2016}F. Lindsten, P. Bunch, S. S\"{a}rkk\"{a}, T. B.
Sch\"{o}n and S. J. Godsill, \textquotedblleft Rao-Blackwellized Particle
Smoothers for Conditionally Linear Gaussian Models\textquotedblright,
\emph{IEEE J. Sel. Topics in Sig. Proc.}, vol. 10, no. 2, pp. 353-365, March 2016.

\bibitem {Loeliger_2007}H.-A. Loeliger, J. Dauwels, Junli Hu, S. Korl, Li
Ping, F. R. Kschischang, ``The Factor Graph Approach to Model-Based Signal
Processing'', \emph{IEEE Proc.}, vol. 95, no. 6, pp. 1295-1322, June 2007.

\bibitem {Kschischang_2001}F. R. Kschischang, B. Frey, and H. Loeliger,
``Factor Graphs and the Sum-Product Algorithm'', \emph{IEEE Trans. Inf.
Theory}, vol. 41, no. 2, pp. 498-519, Feb. 2001.

\bibitem {Berrou_1996}C. Berrou and A. Glavieux, \textquotedblleft Near
Optimum Error Correcting Coding and Decoding: Turbo-Codes\textquotedblright,
\emph{IEEE Trans. Commun.}, vol. 44, no. 10, pp. 1261 - 1271, Oct. 1996.

\bibitem {Benedetto_1998}S. Benedetto, D. Divsalar, G. Montorsi and F.
Pollara, ``Serial Concatenation of Interleaved Codes: Performance Analysis,
Design, and Iterative Decoding'', \emph{IEEE Trans. Inf. Theory}, vol. 44, no.
3, pp. 909-926, May 1998.

\bibitem {Koetter_2004}R. Koetter, A. C. Singer and M. T\"{u}chler, ``Turbo
Equalization'', \emph{IEEE Sig. Proc. Mag.}, vol. 21, no. 1, pp. 67-80, Jan. 2004.

\bibitem {Hagenauer_1997}J. Hagenauer, ``The Turbo Principle: Tutorial
Introduction \& State of the Art\textquotedblright, \emph{Proc. Int. Symp.
Turbo Codes \& Related Topics}, Brest, France, Sep. 1997, pp. 1-11.

\bibitem {Worthen_2001}A. P. Worthen and W. E. Stark, ``Unified Design of
Iterative Receivers using Factor Graphs'', \emph{IEEE Trans. Inf. Theory},
vol. 47, no. 2, pp. 843\textendash849, Feb. 2001.

\bibitem {Kschischang_1998}F. R. Kschischang and B. J. Frey, ``Iterative
Decoding of Compound Codes by Probability Propagation in Graphical Models'',
\emph{IEEE J. Sel. Areas Commun.}, vol. 16, no. 2, pp. 219-230, Feb. 1998.

\bibitem {Dauwels_2006}J. Dauwels, S. Korl and H.-A. Loeliger,
\textquotedblleft Particle Methods as Message Passing\textquotedblright,
\emph{Proc. 2006 IEEE Int. Symp. on Inf. Theory}, pp. 2052-2056, 9-14 July 2006.

\bibitem {Minka_2001}T. P. Minka, \textquotedblleft Expectation Propagation
for Approximate Bayesian Inference\textquotedblright, \emph{Proc. 17th Annual
Conf. Uncertainty in Artif. Intell.}, pp. 362-369, Seattle, WA (USA), Aug. 2001.

\bibitem {Zoeter_2006}O. Zoeter and T. Heskes, ``Deterministic Approximate
Inference Techniques for Conditionally Gaussian State Space Models'',
\emph{Statistics and Computing}, vol. 16, no. 3, pp. 279-292, Sep. 2006.

\bibitem {Smidl_2005}V. Sm\'{\i}dl and A. Quinn, \textbf{The Variational Bayes
Method in Signal Processing}, Berlin, Germany, Springer, 2005.

\bibitem {Bishop}C. M. Bishop, \textbf{Pattern Recognition and Machine
Learning}, Springer, 2006.

\bibitem {Runnalls_2007}A. R. Runnalls, \textquotedblleft Kullback-Leibler
Approach to Gaussian Mixture Reduction\textquotedblright, \emph{IEEE Trans. on
Aer. and Elec. Syst.}, vol. 43, no. 3, pp. 989-999, July 2007.

\bibitem {Li_2015}T. Li, M. Bolic, P. Djuric, \textquotedblleft Resampling
Methods for Particle Filtering: Classification, Implementation, and
Strategies\textquotedblright, \emph{IEEE Sig. Proc. Mag.}, vol. 32, no. 3,
pp.70-86, May 2015.

\bibitem {Vitetta}G. M. Vitetta, D. P. Taylor, G. Colavolpe, F. Pancaldi and
P. A. Martin, \textbf{Wireless Communications: Algorithmic Techniques}, John
Wiley \& Sons, 2013.

\bibitem {Hagenauer_1996}J. Hagenauer, E. Offer and L. Papke,
\textquotedblleft Iterative decoding of binary block and convolutional
codes\textquotedblright, \emph{IEEE Trans. Inf. Theory}, vol. 42, no. 2, pp.
429-445, Mar 1996.

\bibitem {Schon_2010}T. Sch\"{o}n, \textquotedblleft Example Used in
Exemplifying the Marginalized (Rao-Blackwellized) Particle Filter", Nov. 2010
(available at http://users.isy.liu.se/en/rt/schon/Code/RBPF/Document/MPFexample.pdf).
\end{thebibliography}
\end{document}